\documentclass{article}

\usepackage[T1]{fontenc}
\usepackage[utf8]{inputenc}
\usepackage[bookmarks=true, bookmarksopen=true, colorlinks, citecolor=blue, anchorcolor = blue, linkcolor = blue, urlcolor=blue]{hyperref}
\usepackage[utf8]{inputenc}
\usepackage[T1]{fontenc}
\usepackage[english]{babel}
\usepackage{lmodern}
\usepackage{array}
\usepackage{amsmath}
\usepackage{amssymb}
\usepackage{amsfonts}
\usepackage{amsthm}
\usepackage{mathrsfs}
\usepackage{enumitem}
\usepackage{pgf}
\usepackage{bbm}
\usepackage{fullpage}
\usepackage{float}
\usepackage{subfig}
\usepackage[ruled, noend]{algorithm2e}

\title{On determinantal point processes with nonsymmetric kernels}

\author{%
  Arnaud~Poinas\footnote{Laboratoire de Mathematiques et Applications, University of Poitiers, France.}}

\numberwithin{equation}{section}

\renewcommand{\leq}{\leqslant}
\renewcommand{\geq}{\geqslant}
\renewcommand{\tilde}{\widetilde}
\renewcommand{\epsilon}{\varepsilon}

\newcommand{\bC}{\mathbb{C}}

\newcommand{\bE}{\mathbb{E}}

\newcommand{\bK}{\mathbb{K}}

\newcommand{\bN}{\mathbb{N}}

\newcommand{\bP}{\mathbb{P}}

\newcommand{\bR}{\mathbb{R}}

\newcommand{\cM}{\mathcal{M}}

\newcommand{\cP}{\mathcal{P}}

\newcommand{\cS}{\mathcal{S}}

\newcommand{\cara}[1]{\mathbbm{1}_{#1}}

\newcommand{\defeq}{:=}

\newcommand{\Pf}{\textnormal{Pf}}

\newcommand{\DPP}{\textnormal{DPP}}

\theoremstyle{plain}
\newtheorem{proposition}{Proposition}[section]

\newtheorem{definition}[proposition]{Definition}
\newtheorem{defin/prop}[proposition]{Definition/Proposition}
\newtheorem{theorem}[proposition]{Theorem}
\newtheorem{corollary}[proposition]{Corollary}

\begin{document}

\maketitle

\begin{abstract}
Determinantal point processes (DPPs for short) are a common class of repulsive point processes. They have found numerous applications such as the modeling of spatial point pattern datasets with repulsion between close points. In the case of DPPs on finite sets they are defined through a matrix, called the DPP kernel, which is usually assumed to be symmetric. While there are a few known examples of DPPs with nonsymmetric kernels, not much is known on how the removal of the symmetry assumption affects the usual DPP properties. In this paper, we demonstrate how to adapt the results on $P_0$ matrices to the DPP setting in order to get necessary and sufficient conditions for the well-definedness of DPPs with nonsymmetric kernels. We also generalize the common properties of DPPs with symmetric kernel to this broader setting. We then finish by showing how these results can be used to construct attractive couplings of repulsive DPPs.
\end{abstract}

\bigskip

In $1975$, Odile Macchi \cite{Macchi} introduced determinantal point processes (or DPPs for short) to model fermion particles. In this context, DPPs are seen as a random locally finite configuration of points in a continuous space, usually $\bR^d$. Since then, DPPs have seen a wide range of applications not only in physics but also in random matrix theory \cite{Soshnikov} to model the distributions of eigenvalues of some classes of random matrices, in statistics \cite{Lavancier} to model datasets of repulsive point patterns or to generate quadrature points for numerical integration \cite{Bardenet20}, just to name a few examples. More recently, there has been a growing interest in studying a finite discretized version of DPPs defined as a random subset of a finite collection of objects \cite{Taskar, Lyons, ShiraiII}. One of the main application is in machine learning \cite{Taskar} to select a diverse subset of a large dataset. 

\medskip

In almost all cases, determinantal point process are defined through a symmetric kernel. This symmetry is an important assumption needed for a lot of results (well-definedness, simulation, distribution of the number of points, ...), the main one being that DPPs are repulsive point process. In the continuous case it means that the observation of a point at a given location decreases the likelihood of seeing another point nearby. In the discrete case it means that the selection of a given object decreases the likelihood of selecting another object with similar characteristics. Although, the symmetry of the kernel is not always a necessary assumption as some examples of DPPs with nonsymmetric kernels are known \cite{Borodin09, Brunel18, Lyons, Soshnikov}. In the context of machine learning, recent works were especially interested in using nonsymmetric kernels in order to inject some attraction inside the determinantal distribution \cite{Bardenet, Brunel18, Gartrell, Gartrell21, Han, Han22} for applications like recommended systems \cite{Gartrell, Gartrell21}. In particular, \cite{Brunel18, Gartrell, Gartrell21} focus on the inference of nonsymmetric DPPs, \cite{Han, Han22} look at their simulation and \cite{Bardenet} shows concentration inequalities for linear statistics. Nevertheless, these works left a lot of questions open about how the general properties of DPPs are altered in the case of nonsymmetric kernels.

\medskip

% and there has been a few papers working with them \cite{Gartrell, Han}. Nevertheless, these works only focus on some specific questions related to these DPPs but does not establish their general properties. Our goal in this paper is therefore to investigate how the usual results on DPPs are modified when considering a nonsymmetric kernel. 
Our goal in this paper is to extend some of the general results for discrete DPPs with symmetric kernel to the case of generic kernels. We focus exclusively on DPPs on finite sets and the extension of the results thereafter to infinite sets and to the continuous case is left for future works. We begin in Section \ref{sec: Def} by recalling the basic definitions of DPPs over a finite set and some of their main properties that does not require the kernel symmetry to be satisfied. As pointed out in \cite{Brunel18, Gartrell, Han}, discrete DPPs with generic kernels are closely related to the theory of $P_0$-matrices \cite{Johnson20}. We show in Section \ref{sec: Properties} how to translate some of the common results on $P_0$-matrices to the DPP setting. This way, we establish necessary and sufficient conditions on the well-definedness of DPPs with generic kernels as well as various properties on the eigenvalues of these kernels. In particular, theses results points to a natural way of writing a correlation kernel $K$ of a DPP as $K=\frac{1}{2}(I_n - M)$. We focus in Section \ref{sec: prop_DPPs} on the property of this matrix $M$. We show that its principal minors give information about the parity of the number of points of the DPPs. We then prove that, using this matrix $M$, it is possible to generalize the property that all DPPs with symmetric kernels can be writing as a mixture of DPPs with a projection kernel. This is one of the most useful properties of DPPs, used to construct the HKPV algorithm \cite[Algorithm 4.4.2]{Hough} and to characterize DPPs in the continuous case. Finally, we give an application of DPPs with nonsymmetric kernels in Section \ref{sec: Simulation} as a tool to construct and simulate attractive couplings of repulsive DPPs. This can be used as a model for dynamic random sets \cite{Gourieroux} or for marked spatial data with repulsion between points of the same mark and attraction between points of different marks with some examples given at the beginning of the section.

\section{Definitions and first properties of DPPs}\label{sec: Def}

\subsection{Notations and definitions}

For any integer $n\in\bN\backslash\{0\}$, we write $[n]$ for the set $\{1,\cdots, n\}$, $\mathcal{P}([n])$ for the power set of $[n]$, $\mathcal{M}_n(\bR)$ for the set of $n\times n$ matrices with real entries and $\mathcal{S}_n(\bR)$ (resp. $\mathcal{S}^+_n(\bR)$) for the set of symmetric (resp. symmetric positive semi-definite) matrices with real entries. For any finite set $S$, we write $|S|$ for its cardinal. Given any matrix $M\in\mathcal{M}_n(\bR)$ we write $\|M\|_2$ for its spectral norm. For any $S,T\subset[n]$ we define $M_{S,T}$ as the submatrix of $M$ with rows indexed by $S$ and columns indexed by $T$. When $S=T$ we write $M_S$ instead of $M_{S, S}$ in order to simplify the notations. We denote by $I_n$ the $n\times n$ identity matrix and $1_n$ (resp. $0_n$) the vector of $\bR^n$ uniquely composed of $1$s (resp. $0$s). More generally, for any $S\subset [n]$ we write $\cara{S}$ for the vector whose elements indexed by $S$ are equal to $1$ and the others are equal to $0$. For any vector $x=(x_1,\cdots, x_n)^T\in\bR^n$ we write $D(x)$ for the $n\times n$ diagonal matrix whose diagonal elements are $x_1,\cdots, x_n$. Finally, we write $X\sim b(p)$ to indicate that $X$ is a random variable with a Bernoulli distribution with parameter $p\in [0, 1]$.

\medskip

\noindent We begin by recalling the general definition of determinantal measures and $L$-ensemble measures.

\begin{definition}
Let $n\in\bN\backslash\{0\}$. 
\begin{itemize}
\item Let $L\in\mathcal{M}_n(\bR)$ such that $I_n+L$ is invertible. The \textbf{L-ensemble measure} $\mu$ on $\mathcal{P}([n])$ with kernel $L$ is defined by
\[\forall S\subset [n],~\mu(\{S\})=\frac{\det(L_S)}{\det(I_{n}+L)},\]
with the convention $\det(L_\emptyset)=1$.
\item Let $K\in\mathcal{M}_n(\bR)$. The \textbf{determinantal measure} $\mu$ on $\mathcal{P}([n])$ with kernel $K$ is defined by
\[\forall S\subset [n],~\mu(\{X\in\mathcal{P}([n])~\mbox{s.t.}~S\subset X\})=\det(K_S).\]
\end{itemize}
\end{definition}

\noindent As shown in \cite{Taskar}, these two definitions are almost equivalent since, when $I_n-K$ is invertible, a determinantal measure with kernel $K$ is an $L$-ensemble measure with kernel $L=K(I_n-K)^{-1}$. Reciprocally, an $L$-ensemble measure with kernel $L$ is always a determinantal measure with kernel $K=L(I_n+L)^{-1}$. These two measures are always well defined and with total mass equal to $1$ as a consequence of identity \cite[1.2.P20]{Horn}:
\begin{equation}\label{eq: Sum_minors}
\forall M\in\mathcal{M}_n(\bR),~\sum_{S\subset [n]}\det(M_S) = \det(I_n+M).
\end{equation}
Therefore, these measures are probability measures if and only if $\mu(\{S\})\geq 0$ for all $S\subset [n]$. When well-defined, the associated probability distributions are called \textbf{$L$-ensemble} and \textbf{determinantal point processes} and we write them $\DPP_L(L)$ and $\DPP(K)$ respectively. We also write $X\sim \DPP(K)$ (resp. $X\sim \DPP_L(L)$) for a determinantal point process (resp. $L$-ensemble) with kernel $K$ (resp. $L$). When $K$ and $L$ are symmetric matrices then the condition needed for the associated probability measure to be well-defined is already known (see \cite{Taskar} for example).

\begin{proposition}
Let $K,L\in\mathcal{S}_n(\bR)$. The determinantal measure with kernel $K$ is a probability measure if and only if the eigenvalues of $K$ are all in $[0, 1]$. The $L$-ensemble measure with kernel $L$ is a probability measure if and only if $L\in\mathcal{S}_n^+(\bR)$.
\end{proposition}

\noindent For a generic kernel $L$, it has been pointed out before (see \cite{Brunel18, Han} for example) that its associated $L$-ensemble measure is a probability distribution if and only if its principal minors are non-negative. Matrices satisfying this property are called $P_0$-matrices in the linear algebra literature.

\begin{definition}
A matrix $L\in\mathcal{M}_n(\bR)$ is said to be a $P_0$-matrix if
\[\forall S\subset [n],~~\det(L_S)\geq 0.\]
\end{definition}

\noindent We refer to \cite[Chapter 4]{Johnson20} for a review of the general properties of $P_0$-matrices. We give a quick proof of the characterization of $L$-ensemble kernels for the sake of completion.

\begin{proposition}
The $L$-ensemble measure with kernel $L\in\mathcal{M}_n(\bR)$ is a probability measure if and only if $L$ is a $P_0$ matrix.
\end{proposition}
\begin{proof}
We already mentioned that the $L$-ensemble measure with kernel $L$ is a probability measure if and only if
\begin{equation}\label{eq: L_prob}
\forall S\subset [n],~~\mu(\{S\})=\frac{\det(L_S)}{\det(I_n+L)}\geq 0.
\end{equation}
Taking $S=\emptyset$ shows that we need $\det(I_n+L)$ to be positive hence \eqref{eq: L_prob} is equivalent to $\det(L_S)\geq 0$ for all $S\subset [n]$ and $\det(I_n+L)\geq 0$. Using \eqref{eq: Sum_minors} we can then conclude that \eqref{eq: L_prob} is simply equivalent to $L$ having its principal minors be non-negative.
\end{proof}

\noindent We note that, since the main focus in machine learning applications is on $L$-ensembles, some special cases of $P_0$-matrices have already been used to construct nonsymmetric $L$ kernels. For example, \cite{Gartrell} consider $L$-ensembles where 
\begin{equation} \label{eq: L_PSD}
L + L^T\in \cS_n^+(\bR)
\end{equation}
and \cite{Gartrell21, Han, Han22} consider $L$-ensembles of the form $L=VV^T + B(D-D^T)B^T$. On the opposite, there is no already well-known matrix family corresponding to the set of kernels $K$ whose associated determinantal measure is a probability distribution. In the rest of a paper we call such matrices a \textbf{DPP kernel} and we focus mainly on them instead of $L$-ensembles.

\subsection{Properties of DPPs not needing the kernel symmetry}

We recall some of the well-known standard properties of DPPs that can be proved without using the assumption that their kernel is symmetrical.

\begin{proposition}\label{prop: Already_proved}
Let $K\in\mathcal{M}_n(\bR)$ be a DPP kernel with complex eigenvalues (taken with multiplicities) $\lambda_1,\cdots,\lambda_n\in\bC$ and define $X\sim \DPP(K)$. Then,
\begin{enumerate}[label=(\alph*)]
\item For any $S\subset [n]$, $X\cap S\sim \DPP(K_S)$. In particular, $K_S$ is a DPP kernel. \label{subprop: Marginal}
\item Let $p\in[0, 1]$ and let $Y$ be the $p$-thinning of $X$ (~$Y$ is  obtained by removing each point of $X$ independently with probability $1-p$) then $Y\sim \DPP(pK)$ and thus $pK$ is a DPP kernel. \label{subprop: Thinning}
\item $X^c\sim \DPP(I_n - K)$. In particular, $I_n - K$ is a DPP kernel. \label{subprop: Complement}
\item If $K$ is a block diagonal matrix with diagonal blocks $K_1,\cdots, K_l$ then $K$ being a DPP kernel is equivalent to each $K_i$ being a DPP kernel. \label{subprop: Block_diag}
\item 
\[\bE\left[\binom{|X|}{k}\right]=e_k(\lambda_1,\cdots,\lambda_n),\]
where $e_k$ is the $k$-th elementary symmetric polynomial. \label{subprop: Moments_nbpoint}
\item For any $x=(x_1,\cdots, x_n)^T\in\bR^n$,
\[\bE\left[\prod_{i\in X} x_i\right]=\det(I_n-D(1_n-x)K).\]
\label{subprop: Laplace}
\item For any $i\in [n]$ such that $K_{i,i}\neq 0$, the distribution of $X\backslash \{i\}$ conditionally to $i\in X$ is a DPP with kernel
\begin{equation}\label{eq: conditional_K}
\tilde K = K-\frac{1}{K_{i,i}}K_{[n], i}K_{i, [n]}.
\end{equation}
In particular, $\tilde K$ is a DPP kernel. \label{subprop: Cond_point}
\end{enumerate}
\end{proposition}
\begin{proof}
\ref{subprop: Marginal} is an immediate consequence of the definition of DPPs. \ref{subprop: Thinning} is the direct result of $\det((pK)_S)=p^{|S|}\det(K_S)$. \ref{subprop: Complement} comes from identity \eqref{eq: Sum_minors} applied to $K_S$:
\[\det((I_n-K)_S)=\sum_{T\subset S}(-1)^{|T|}\det(K_T)\]
combined with the inclusion-exclusion principle. \ref{subprop: Block_diag} comes directly from the fact that the determinant of a block diagonal matrix is the product of the determinants of its blocks. For \ref{subprop: Moments_nbpoint}, note that $\binom{|X|}{k}$ is the number of subsets of $[n]$ with size $k$ thus 
\[\bE\left[\binom{|X|}{k}\right]=\sum_{\substack{S\subset [n]\\ |S|=k}}\mathbb{P}(S\subset X)=\sum_{\substack{S\subset [n]\\ |S|=k}}\det(K_S).\]
It is well-known \cite[Theorem 1.2.16]{Horn} that the sum of principal minors of size $k$ of any matrix $K$ is equal to $e_k(\lambda_1,\cdots,\lambda_n)$ hence the result. \ref{subprop: Laplace} is a particular case of the Laplace transform of DPPs shown in \cite{Shirai} and can be proved using \eqref{eq: Sum_minors} and the decomposition
\[\prod_{i\in X} x_i=\prod_{i\in X} (1-(1-x_i))=\sum_{S\subset [n]}(-1)^{|S|}\left(\prod_{i\in S}(1-x_i)\right)\cara{S\subset X}.\]
Finally, \ref{subprop: Cond_point} is a direct application of the Schur complement (see \cite[Theorem 6.5]{Shirai} for example).
\end{proof}

\subsection{The particle-hole involution and the principal pivot transform}

One of the earliest known way to create DPPs with nonsymmetric kernels is to use what is often called the \textbf{particle-hole transformation}. For any subset $S$ of $[n]$, the particle-hole transformation with respect to a set $S\subset [n]$ is the involution $X \mapsto (X\cap S^c)\cup(X^c\cap S)$ that switches the states of the points in and out of $X\cap S$. It is shown in \cite{Borodin99} that determinantal distributions are stable by this operation.

\begin{proposition}[\cite{Borodin99}]\label{prop: Particle_Hole}
Let $K\in\mathcal{M}_n(\bR)$ be a DPP kernel and let $X\sim \DPP(K)$. For a given set $S\subset [n]$, we define
\[\tilde X = (X\cap S^c) \cup (X^c\cap S).\]
Then, $\tilde X\sim \DPP(\tilde K)$ with
\begin{equation}\label{eq: PH_K}
\tilde K=D(\cara{S})(I_n-K)+D(\cara{S^c})K.
\end{equation}
After some permutations of the rows and columns of $K$ and $\tilde K$ we can write
\begin{equation}\label{eq: Particle_Hole_rowcol_permut}
K=\begin{pmatrix} K_S & K_{S, S^c} \\
K_{S^c, S} & K_{S^c}
\end{pmatrix}~\Rightarrow~\tilde K=\begin{pmatrix} I_{|S|}-K_S & -K_{S, S^c} \\
K_{S^c, S} & K_{S^c}
\end{pmatrix}.
\end{equation}
\end{proposition}

\noindent As a direct consequence we get an expression of all probabilities $\mathbb{P}(X=S)$ using only the kernel $K$:

\begin{corollary}[Identity (147) in \cite{Taskar} and Identity (2.1) in \cite{ShiraiII}]\label{corr: Alternate_LL}
\[\mathbb{P}(X=S^c)=\mathbb{P}([n]\subset \tilde X)=\det(D(\cara{S})(I_n-K)+D(\cara{S^c})K)=\det\begin{pmatrix} I_{|S^c|}-K_{S^c} & -K_{S^c, S}\\
K_{S, S^c} & K_S
\end{pmatrix}.\]
\end{corollary}

\noindent We now give a generalization of Proposition \ref{prop: Particle_Hole} where, instead of switching the state of a given set of points, we choose whether we switch the state of each point in $[n]$ or not independently with its own probability.

\begin{proposition}\label{prop: Switching_the_values}
Let $K\in\mathcal{M}_n(\bR)$ be a DPP kernel and let $X\sim \DPP(K)$. We define $p=(p_1,\cdots,p_n)\in[0,1]^n$ and we consider some independent Bernoulli random variables $B_i\sim b(p_i)$ also independent from $X$. We construct $\tilde X$ as
\[\tilde X=\{i\in [n]~\mbox{s.t.}~i\in X~\mbox{and}~B_i=0~\mbox{or}~i\notin X~\mbox{and}~B_i=1\}.\]
Then,
\[\tilde X\sim \DPP(D(p)(I_n-K)+D(1_n-p)K)\]
\end{proposition}
\begin{proof}
For any $S\subset [n]$ we write
\begin{align*}
\mathbb{P}(S\subset \tilde X)&=\sum_{T\subset S}\mathbb{P}(S\subset \tilde X|B_i=1~\mbox{when}~i\in T~\mbox{and}~B_i=0~\mbox{when}~i\in T\backslash S)\prod_{i\in T}p_i\prod_{i\in S\backslash T}(1-p_i)\\
&=\sum_{T\subset S}\mathbb{P}(X\cap S=S\backslash T)\prod_{i\in T}p_i\prod_{i\in S\backslash T}(1-p_i).
\end{align*}
Using Corollary \ref{corr: Alternate_LL} we get
\begin{align*}
\mathbb{P}(S\subset \tilde X)&=\sum_{T\subset S}\det\begin{pmatrix} I_{|T|}-K_{T} & -K_{T, S\backslash T}\\
K_{S\backslash T, T} & K_{S\backslash T}
\end{pmatrix}\prod_{i\in T}p_i\prod_{i\in S\backslash T}(1-p_i)\\
&=\sum_{T\subset S}\det\left(\begin{pmatrix} D(p)_T & 0 \\
0 & D(1_n-p)_{S\backslash T}
\end{pmatrix}\begin{pmatrix} I_{|T|}-K_{T} & -K_{T, S\backslash T}\\
K_{S\backslash T, T} & K_{S\backslash T}
\end{pmatrix}\right)
\end{align*}
Since $D(p)$ and $D(1_n-p)$ are diagonal then the matrix inside the determinant in the last identity has its rows indexed by the elements of $T$ corresponding to the rows of $(D(p)(I_n-K))_S$ and its rows indexed by the elements of $S\backslash T$ corresponding to the rows of $(D(1_n-p)K)_S$. Therefore, by the multilinearity of the determinant with respect to its rows we get that
\[\mathbb{P}(S\subset \tilde X)=\det\big((D(p)(I_n-K)+D(1_n-p)K)_S\big),\]
concluding the proof.
\end{proof}

\noindent Note that we recover Proposition \ref{prop: Particle_Hole} by taking $p=\cara{S}$ for some $S\subset [n]$. We now show that the particle hole involution is closely linked to a common tool appearing in the theory of $P_0$ matrices called the \textbf{principal pivot transform} \cite{Tsatsomeros}.

\begin{definition}\label{def: PPT}
Let $M\in\mathcal{M}_n(\bR)$ and $S\subset [n]$ such that $M_S$ is invertible. With the right permutation of rows and columns we can write $M$ as $\begin{pmatrix} M_S & M_{S, S^c} \\
M_{S^c, S} & M_{S^c} \end{pmatrix}.$ The \textbf{principal pivot transform} of $M$ relative to $S$ is then defined as
\[ppt(M,S)\defeq\begin{pmatrix} M_S^{-1} & -M_S^{-1}M_{S, S^c} \\
M_{S^c, S}M_S^{-1} & M_{S^c} - M_{S^c, S}M_S^{-1}M_{S, S^c} \end{pmatrix}.\]
\end{definition}

\noindent The definition of the principal pivot transform is often better understood through the following result.

\begin{proposition}\label{prop: PPT_def_alternate}
Let $M\in\mathcal{M}_n(\bR)$ and $S\subset [n]$ such that $M_S$ is invertible and $x,y\in \bR^n$. We use the same permutation as in Definition \ref{def: PPT} and write $x=\begin{pmatrix} x_S \\ x_{S^c} \end{pmatrix}$ and $y=\begin{pmatrix} y_S \\ y_{S^c} \end{pmatrix}$. Then,
\[\begin{pmatrix} y_S \\ y_{S^c} \end{pmatrix} = M\begin{pmatrix} x_S \\ x_{S^c} \end{pmatrix}~\Leftrightarrow~\begin{pmatrix} x_S \\ y_{S^c} \end{pmatrix}=ppt(M,S)\begin{pmatrix} y_S \\ x_{S^c} \end{pmatrix}.\]
\end{proposition}

\noindent We now show the link between the principal pivot transform and DPPs through the following result.

\begin{proposition}\label{prop: linear_PPT}
Let $L\in\mathcal{M}_n(\bR)$ be a $P_0$ matrix such that $L_S$ is invertible and let $X\sim \DPP_L(L)$. For a given set $S\subset [n]$ we define the particle-hole transformation of $X$ as
\[\tilde X = (X\cap S^c) \cup (X^c\cap S).\] 
If $I_n-\tilde K$ is invertible then $\tilde X\sim \DPP_L(ppt(L, S))$.
\end{proposition}
\begin{proof}
Let $K=L(I_n+L)^{-1}$ be the DPP kernel of $X$ and $\tilde K$ be the DPP kernel of $\tilde X$. After some row and column permutations we write $K$ and $\tilde K$ as in identity \eqref{eq: Particle_Hole_rowcol_permut}. Since $I_n-\tilde K$ is assumed to be invertible, we can define $\tilde L=(I_n-\tilde K)^{-1}-I_n$. We then choose $x,y\in\bR^n$ indexed as in Proposition \ref{prop: PPT_def_alternate} such that $y=\tilde L x$. Then,
\[y=\tilde L x~\Leftrightarrow~(x+y)=(I_n+\tilde L)x~\Leftrightarrow~(I_n-\tilde K)(x+y)=x~\Leftrightarrow~y=\tilde K(x+y).\]
This expression can be extended into
\begin{align*}
& \begin{pmatrix} y_S \\ y_{S^c} \end{pmatrix}=\begin{pmatrix} x_S + y_S -K_S(x_S+y_S) -K_{S, S^c}(x_{S^c}+y_{S^c}) \\
K_{S^c, S}(x_S+y_S) + K_{S^c}(x_{S^c}+y_{S^c})
\end{pmatrix}\\
\Leftrightarrow~ & \begin{pmatrix} -x_S \\ y_{S^c} \end{pmatrix}=\begin{pmatrix} -K_S(x_S+y_S) -K_{S, S^c}(x_{S^c}+y_{S^c}) \\
K_{S^c, S}(x_S+y_S) + K_{S^c}(x_{S^c}+y_{S^c})
\end{pmatrix}\\
\Leftrightarrow~ & \begin{pmatrix} x_S \\ y_{S^c} \end{pmatrix}=K(x+y)\\
\Leftrightarrow~ & \begin{pmatrix} y_S \\ x_{S^c} \end{pmatrix}=(I_n-K)(x+y)\\
\Leftrightarrow~ & (I_n+L)\begin{pmatrix} y_S \\ x_{S^c} \end{pmatrix}=x+y\\
\Leftrightarrow~ & \begin{pmatrix} x_S \\ y_{S^c} \end{pmatrix} = L\begin{pmatrix} y_S \\ x_{S^c} \end{pmatrix}.
\end{align*}
We conclude that $\tilde L = ppt(L, S)$ with Proposition \ref{prop: PPT_def_alternate}.
\end{proof}

\noindent The result that determinantal measures are stable by particle-hole transformations \cite{Borodin99} is thus almost equivalent to the result that $P_0$ matrices are stable by principal pivot transforms \cite[Theorem 4.8.4]{Johnson20}. Interestingly, both results were proven almost at the same time but using completely different methods. This is a nice illustration of the close link there is between the theory of $P_0$ matrices and the theory of DPPs.

\section{Translating the theory of \texorpdfstring{$P_0$}{P0}-matrices into the DPP setting}\label{sec: Properties}

\subsection{Characterization of DPP kernels}

The first hurdle in working with nonsymmetric kernels is the issue of verifying whether a given matrix is a DPP kernel or not. In the case of $L$-ensemble, there is already a lot of literature on the problem of testing if a matrix is $P_0$ or not and we refer to \cite[Section 4.6]{Johnson20} for an overview. We focus instead on adapting the results on $P_0$ matrices to characterize whether a given matrix $K$ is a DPP kernel or not. We first recall that it was proved in \cite{Coxson} that the problem of checking whether a given matrix is a $P$ matrix (matrix with positive principal minors) is co-NP-complete. The same can be expected for checking if a given matrix is a DPP kernel since we can go from $K$ to $L$ in polynomial time. We thus can't expect any characterization of DPP kernels that is as easy to check as just computing the eigenvalues of $K$ like in the symmetrical case. Instead, we give the following characterizations based on \cite[Theorem 4.3.4 and Theorem 4.3.9]{Johnson20}.

\begin{theorem}\label{theo: Caracterisation_DPPs}
Let $K\in\mathcal{M}_n(\bR)$. $K$ is a DPP kernel if and only if one of the following equivalent properties is satisfied:
\begin{enumerate}[label=(\alph*)]
\item \label{cara1} \begin{equation}\label{eq: cara1}
\forall p\in\{0, 1\}^n,~\det(D(p)(I_n-K)+D(1_n-p)K)\geq 0.
\end{equation}
\item \label{cara2} \begin{equation}\label{eq: cara2}
\forall p\in (0, 1)^n,~\det(D(p)(I_n-K)+D(1_n-p)K)> 0.
\end{equation}
\item For all non-zero $x\in\bR^n$ there exists $i\in [n]$ such that 
\[x_i(Kx)_i\geq 0,~|(Kx)_i|\leq|x_i|~\mbox{and}~x_i \neq 0.\] \label{cara3}
\end{enumerate}
\end{theorem}
\begin{proof}
The proof of \ref{cara3} is directly adapted from the method used in \cite[Theorem 4.3.4]{Johnson20} while the proofs of \ref{cara1} and \ref{cara2} are done differently using the interpretation of \eqref{eq: cara1} and \eqref{eq: cara2} in the setting of DPPs.
\begin{enumerate}[label=(\alph*)]
\item Direct consequence of Corollary \ref{corr: Alternate_LL}.
\item We first remark that $D(p)(I_n-K)+D(1_n-p)K$ is invertible for all $p\in (0, 1)^n$ if and only if $\det(D(p)(I_n-K)+D(1_n-p)K)>0$ for all $p\in (0, 1)^n$. This is a consequence of the continuity of the determinant and the fact that if $p=\frac{1}{2}1_n$ then
\[\det(D(p)(I_n-K)+D(1_n-p)K) = \det\left(\frac{1}{2}I_n\right)>0.\]
The continuity of the determinant also shows that $\ref{cara2}\Rightarrow\ref{cara1}\Rightarrow K$ is a DPP kernel.
Now, let $K$ be a DPP kernel, choose $p\in (0, 1)^n$, define $X\sim \DPP(K)$ and construct $\tilde X$ as in Proposition \ref{prop: Switching_the_values}. $\tilde X$ is then a determinantal point process with kernel $D(p)(I_n-K)+D(1_n-p)K$. In particular, we have
\[\mathbb{P}(\tilde X=[n])=\mathbb{P}([n]\subset\tilde X)=\det\big(D(p)(I_n-K)+D(1_n-p)K\big).\]
Now, let $S\subset [n]$ such that $\mathbb{P}(X=S)>0$. Then, by definition of $\tilde X$ we have
\[\mathbb{P}(\tilde X=[n])\geq \mathbb{P}(\tilde X=[n]|X=S)\mathbb{P}(X=S)=\prod_{i\notin S}p_i\prod_{i\in S}(1-p_i)\mathbb{P}(X=S)>0.\]
Therefore, $\det(D(p)(I_n-K)+D(1_n-p)K)>0$ for all $p\in (0, 1)^n$.

\item We assume that there exists $p\in (0, 1)^n$ such that $D(p)(I_n-K)+D(1_n-p)K$ is not invertible. This means that for some non-zero $x\in\bR^n$ we have
\begin{align*}
&D(p)(I_n-K)x+D(1_n-p)Kx = 0\\
\Leftrightarrow~&D(1_n-2p)Kx=-D(p)x\\
\Leftrightarrow~&x=\begin{pmatrix}
2-\frac{1}{p_1} & & 0 \\
& \ddots & \\
& & 2-\frac{1}{p_n}
\end{pmatrix}Kx.
\end{align*}
Since the function $t\mapsto 2-\frac {1}{t}$ is a bijection from $(0, 1)$ to $(-\infty, 1)$ then we can conclude that $K$ is not a DPP kernel if and only if there exists a non-zero $x\in\bR^n$ such that
\[\forall i\in[n],~\exists \mu_i\in (-\infty, 1)~\mbox{s.t.}~x_i=\mu_i(Kx)_i\]
which is equivalent to
\[\forall i\in[n],~x_i(Kx)_i < 0~\mbox{or}~|(Kx)_i| > |x_i|~\mbox{or}~x_i=0.\]
As a consequence, $K$ is a DPP kernel iff for all non-zero $x\in\bR^n$ there exists $i\in [n]$ such that
\begin{equation*}
x_i(Kx)_i \geq 0~,~|(Kx)_i| \leq |x_i|~\mbox{and}~x_i\neq 0.\qedhere
\end{equation*}
\end{enumerate}
\end{proof}

\noindent We note that a different characterization of generic DPP kernels was given in \cite{Brunel} as 
\[K~\textnormal{is a DPP kernel}~\Leftrightarrow~\forall S\subset [n],~(-1)^{|S|}\det(K - D(\cara{S})) \geq 0.\]
This actually corresponds to Theorem \ref{theo: Caracterisation_DPPs} \ref{cara1} since $(-1)^{|S|} = \det(I_n - 2D(\cara{S}))$ hence
\[(-1)^{|S|}\det(K - D(\cara{S})) = \det(D(\cara{S})(I_n-K)+D(1_n-\cara{S})K).\]

\noindent A direct consequence of these characterizations is that the set of DPP kernels is a star-shaped set centered at $\frac 1 2 I_n.$

\begin{proposition}\label{prop: star-shaped}
Let $K\in\mathcal{M}_n(\bR)$ be a DPP kernel. Then, for all $\lambda\in [0, 1]$, $(1-\lambda) K + \lambda\left(\frac 1 2 I_n\right)$ is a DPP kernel.
\end{proposition}
\begin{proof}
Let $\tilde K=(1-\lambda) K + \lambda\left(\frac 1 2 I_n\right)$ and $p\in (0, 1)^n$. By Theorem \ref{theo: Caracterisation_DPPs}~\ref{cara2} we need to show that $D(p)(I_n-\tilde K) + D(1_n-p)\tilde K$ is invertible to conclude that $\tilde K$ is a DPP kernel. We write
\begin{align*}
D(p)(I_n-\tilde K) + D(1_n-p)\tilde K&= D(p)I_n+D(1_n-2p)\tilde K\\
&= D(p)I_n+D(1_n-2p)\left((1-\lambda) K + \frac \lambda 2 I_n\right)\\
&= D\left(p + \frac \lambda 2 1_n - \lambda p\right)I_n+D(1_n-2p-\lambda 1_n + 2\lambda p)K\\
&= D\left(p + \frac \lambda 2 1_n - \lambda p\right)(I_n-K) + D\left(1_n-\left(p + \frac \lambda 2 1_n - \lambda p\right)\right)K
\end{align*}
Now, for all $i\in [n],$
\[\left|p_i + \frac \lambda 2 - \lambda p_i - \frac 1 2\right|=|1-\lambda|\left|\frac 1 2 - p_i\right|< \frac 1 2\]
hence $p_i + \frac \lambda 2 - \lambda p_i\in (0, 1)$ and thus $D(p)(I_n-\tilde K) + D(1_n-p)\tilde K$ is invertible by Theorem \ref{theo: Caracterisation_DPPs}~\ref{cara2}, concluding the proof.
\end{proof}

\subsection{Generating DPP kernels}

It is known that row-diagonally dominant matrices are $P_0$ \cite[Proposition 4.5.1]{Johnson20} and we can get a similar result for DPP kernels.

\begin{proposition}
Let $K\in\mathcal{M}_n(\bR)$ such that $K$ and $I_n-K$ are row diagonally dominant, meaning that the diagonal elements of $K$ are in $[0,1]$ and satisfy
\[\forall i\in [n],~\min(K_{i,i}, 1-K_{i,i})\geq\sum_{j\neq i} |K_{i,j}|.\]
Then $K$ is a DPP kernel.
\end{proposition}
\begin{proof}
Let $p\in \{0,1\}^n$ and $\tilde K=D(p)(I_n-K)+D(1_n-p)K$.
For all distinct $i,j\in[n]$ we have, 
\[\tilde K_{i,j}=p_i(-K_{i,j})+(1-p_i)K_{i,j}.\] 
Then,
\[\sum_{j\neq i}|\tilde K_{i,j}|\leq p_i \sum_{j\neq i}|K_{i,j}|+(1-p_i)\sum_{j\neq i}|K_{i,j}|\leq p_i (1-K_{i,i})+(1-p_i)K_{i,i}=\tilde K_{i,i}.\]
Using the Gershgorin circle theorem we get that any real eigenvalue of $\tilde K$ is non-negative hence $\det(\tilde K)\geq 0$ and by Theorem \ref{theo: Caracterisation_DPPs} \ref{cara1} we can conclude that $K$ is a DPP kernel.
\end{proof}

\noindent It is also known that matrices $L\in\cM_n(\bR)$ satisfying $\langle x, Lx\rangle \geq 0$ for all $x\in\bR^n$ are $P_0$ \cite[Proposition 4.5.2]{Johnson20}. This is equivalent to \eqref{eq: L_PSD} which is the setting used for $L$-ensemble in some previous works like \cite{Gartrell21, Han}. Adapting this result for DPP kernels yields the following proposition.

\begin{proposition}\label{prop: perturb_half_identity}
Let $M\in\mathcal{M}_n(\bR)$ such that $\|M\|_2\leq 1$. Then, $K\defeq \frac{1}{2}(I_n-M)$ is a DPP kernel.
\end{proposition}
\begin{proof}
We can write for all $x\in\bR^n$,
\begin{equation}\label{eq: scalar_product_bound_1}
\langle Kx, (I_n-K)x\rangle = \frac{1}{4}\langle (I_n-M)x, (I_n+M)x\rangle
=\frac{1}{4}\left(\|x\|^2 - \langle x, M^TMx\rangle\right)\geq 0.
\end{equation}
This is a consequence of the largest eigenvalue of $M^TM$ being $\|M\|^2_2=1$ from our assumptions. Now, if $K$ is not a DPP kernel then by Theorem \ref{theo: Caracterisation_DPPs}~\ref{cara2} there exists a non-zero $x\in\bR^n$ and some $p\in (0, 1)^n$ such that
\[D(p)(I_n-K)x+D(1_n-p)Kx=0.\]
In particular, for all $i\in [n]$ we have
\begin{equation}\label{eq: scalar_product_bound_2}
p_i((I_n-K)x)_i+(1-p_i)(Kx)_i=0
\end{equation}
hence $(Kx)_i$ and $((I_n-K)x)_i$ are either both equal to $0$ or of opposite sign and thus $\langle Kx, (I_n-K)x\rangle \leq 0$. Identity \eqref{eq: scalar_product_bound_1} then gives $\langle Kx, (I_n-K)x\rangle = 0$ and by \eqref{eq: scalar_product_bound_2} we must have $Kx=(I_n-K)x=0$ and thus $x\in\textnormal{Ker}(K)\cap\textnormal{Ker}(I_n-K)=\{0_n\}$ which contradicts the assumption that $x\neq 0_n$ and therefore proves the proposition.
\end{proof}

\noindent This result shows that all matrices close to $\frac{1}{2}I_n$ are DDP kernels. Propositions \ref{prop: star-shaped} and \ref{prop: perturb_half_identity} both illustrate that $\frac{1}{2}I_n$ appears as the center of the set of DPP kernels and suggests that it is natural to write DPP kernels as $\frac{1}{2}(I_n-M)$. We later show in section \ref{sec: prop_DPPs} that, when writing DPP kernels this way, this matrix $M$ satisfies a lot of useful properties.

\subsection{Eigenvalues of DPP kernels}

The eigenvalues of DPP kernels are an import tool to study DPPs. We first recall that the well-definedness of DPPs with symmetric kernels is characterized by the eigenvalues being in $[0, 1]$. Moreover, as a consequence of Proposition \ref{prop: Already_proved}~\ref{subprop: Moments_nbpoint}, we know that all moments of the number of points of a DPP only depends on the eigenvalues of its kernel. It means that the eigenvalues of a DPP kernel $K$ (or $L$) fully characterize the distribution of the cardinal of the associated DPP. For these reasons, we focus in this section on studying the behaviour of the eigenvalues of generic DPP kernels. We begin by giving some bound on these eigenvalues using the result of \cite{Kellogg} on $P_0$ matrices we recall below.

\begin{proposition}[\cite{Kellogg}] \label{prop: ev_Lens}
Let $L\in\mathcal{M}_n(\bR)$ be a $P_0$ matrix and $\lambda\in\mathbb{C}^*$ be a non-zero eigenvalue of $L$. Then,
\[|\arg(\lambda)|\leq\pi - \frac{\pi}{n}.\]
In particular, any real eigenvalue of $L$ is in $\bR_+$.
\end{proposition}

\noindent Adapting this result in the determinantal setting gives the following bounds on the eigenvalues of a DPP kernel.

\begin{proposition} \label{prop: ev_DPP}
Let $\lambda$ be an eigenvalue of a DPP kernel $K\in\mathcal{M}_n(\bR)$. Then,
\begin{equation}\label{eq: ev_DPP}
\lambda\in \mathcal{B}_\mathbb{C}\left(\frac{1}{2}+\frac{1}{2\tan\left(\frac \pi n\right)}i,~\frac{1}{2\sin(\frac \pi n)}\right)\cup \mathcal{B}_\mathbb{C}\left(\frac{1}{2}-\frac{1}{2\tan\left(\frac \pi n\right)}i,~\frac{1}{2\sin(\frac \pi n)}\right),
\end{equation}
where we denote by $\mathcal{B}_\mathbb{C}(x,r)$ the complex ball centered in $x$ with radius $r$. In particular, any real eigenvalue of a DPP kernel lies in $[0, 1]$.
\end{proposition}
\begin{proof}
We begin by pointing out that when $n=1$ then \eqref{eq: ev_DPP} corresponds to $\lambda\in [0, 1]$ which is obvious since the only DPP kernels of size $1\times 1$ have values in $[0, 1]$. We now consider that $n\geq 2$ and first assume that $I_n-K$ is invertible. In that case, $\lambda/(1-\lambda)$ is an eigenvalue of $K(I_n-K)^{-1}$ which is a $P_0$ matrix and by Proposition \ref{prop: ev_Lens} we get that either $\lambda=0$ or
\begin{equation}\label{eq: cond_eig_K}
\left|\arg\left(\frac{\lambda}{1-\lambda}\right)\right|\leq\pi - \frac{\pi}{n}.
\end{equation}
If $\lambda\in\bR$ then $\frac{\lambda}{1-\lambda}$ is also purely real and thus it is non-negative by Proposition \ref{prop: ev_Lens} which is only possible when $\lambda\in [0, 1]$. Otherwise, if $\lambda\in\bC\backslash\bR$ then $\bar\lambda$ is also an eigenvalue of $K$ so we first assume that $\Im(\lambda)>0$ and write $\lambda=a+bi$ with $a\in\bR$ and $b>0$. Then,
\[\frac{\lambda}{1-\lambda}=\frac{\lambda-|\lambda|^2}{|1-\lambda|^2}~\Rightarrow~\arg\left(\frac{\lambda}{1-\lambda}\right)=\arg(\lambda-|\lambda|^2)=\arg(a-(a^2+b^2)+bi).\]
Since $b>0$ then
\[\arg(a-(a^2+b^2)+bi)=\left\{\begin{array}{l}
\arctan\left(\frac{b}{a-(a^2+b^2)}\right)~\mbox{if}~a>a^2+b^2.\\
\frac{\pi}{2}~\mbox{if}~a=a^2+b^2.\\
\arctan\left(\frac{b}{a-(a^2+b^2)}\right)+\pi~\mbox{if}~a<a^2+b^2.\\
\end{array}
\right.\]
If $a\geq a^2+b^2$ (equivalent to $\lambda\in \mathcal{B}_\mathbb{C}\left(\frac 1 2, \frac 1 2\right)$) then $|\arg(a-(a^2+b^2)+bi)|\leq \pi/2$ and thus $\left|\arg\left(\frac{\lambda}{1-\lambda}\right)\right|\leq\frac{\pi}{2}\leq \pi - \frac{\pi}{n}$. If $a< a^2+b^2$ then 
\begin{align*}
\left|\arg\left(\frac{\lambda}{1-\lambda}\right)\right|\leq\pi - \frac{\pi}{n}~&\Leftrightarrow~\arctan\left(\frac{b}{a-(a^2+b^2)}\right)\leq -\frac{\pi}{n}\\
&\Leftrightarrow~\frac{b}{a-(a^2+b^2)}\leq -\tan\left(\frac{\pi}{n}\right)\\
&\Leftrightarrow~\frac{b}{\tan\left(\frac{\pi}{n}\right)}\geq (a^2+b^2)-a\\
&\Leftrightarrow~(a^2 - a)+\left(b^2-\frac{b}{\tan\left(\frac{\pi}{n}\right)}\right)\leq 0\\
%&\Leftrightarrow~\frac{1}{4}+\frac{1}{4\tan\left(\frac{\pi}{n}\right)^2}\geq \left(a-\frac 1 2\right)^2 + \left(b - \frac{1}{\tan\left(\frac \pi n\right)}\right)^2\\
&\Leftrightarrow~\left(a-\frac 1 2\right)^2 + \left(b - \frac{1}{2\tan\left(\frac \pi n\right)}\right)^2 \leq \frac{1}{4\sin\left(\frac{\pi}{n}\right)^2}.
\end{align*}
Now, note that $\mathcal{B}_\mathbb{C}(1/2, 1/2)\subset \mathcal{B}_\mathbb{C}\left(\frac{1}{2}+\frac{1}{2\tan\left(\frac \pi n\right)}i,~\frac{1}{2\sin(\frac \pi n)}\right)$ hence
\[\lambda\in\mathcal{B}_\mathbb{C}\left(\frac{1}{2}+\frac{1}{2\tan\left(\frac \pi n\right)}i,~\frac{1}{2\sin(\frac \pi n)}\right).\]
Since $\bar\lambda$ also satisfies \eqref{eq: cond_eig_K} then with the same reasoning we get the desired result. Finally, if $I_n-K$ is not invertible then there exists some $\epsilon>0$ such that $I_n-pK$ is invertible for all $p\in[1-\epsilon, 1)$. 
Since $pK$ is a DPP kernel by Proposition \ref{prop: Already_proved} \ref{subprop: Thinning} then
\[p\lambda\in \mathcal{B}_\mathbb{C}\left(\frac{1}{2}+\frac{1}{2\tan\left(\frac \pi n\right)}i,~\frac{1}{2\sin(\frac \pi n)}\right)\cup \mathcal{B}_\mathbb{C}\left(\frac{1}{2}-\frac{1}{2\tan\left(\frac \pi n\right)}i,~\frac{1}{2\sin(\frac \pi n)}\right)\]
for all $p\in[1-\epsilon, 1)$ and we get the desired result by taking $\epsilon\rightarrow 0$.
\end{proof}

\begin{figure}[htbp]
\centering
\subfloat[Eigenvalues of $P_0$ matrices]{\includegraphics[width=8cm]{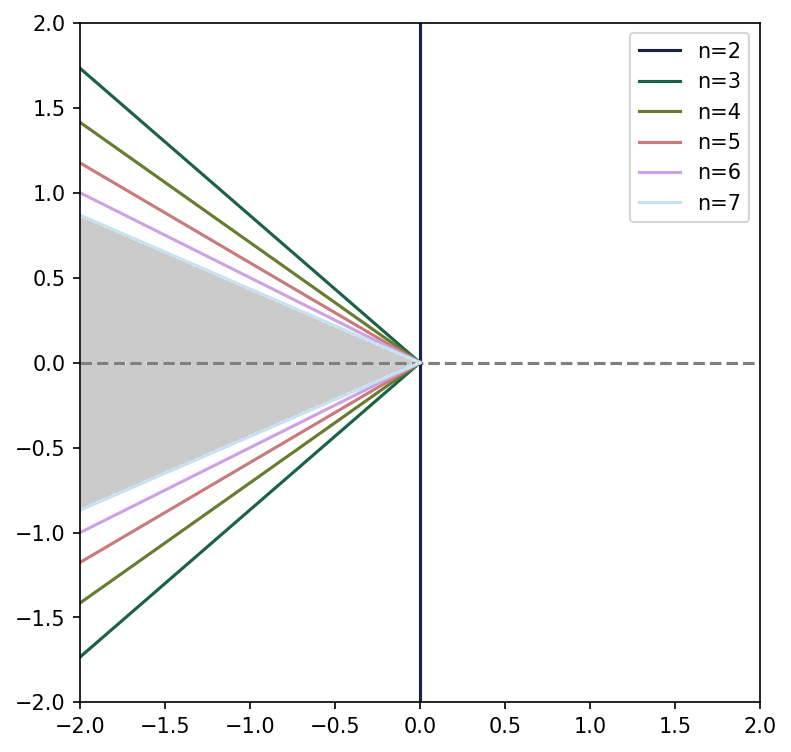}}
\hfill
\subfloat[Eigenvalues of DPP kernels.]{\includegraphics[width=6cm]{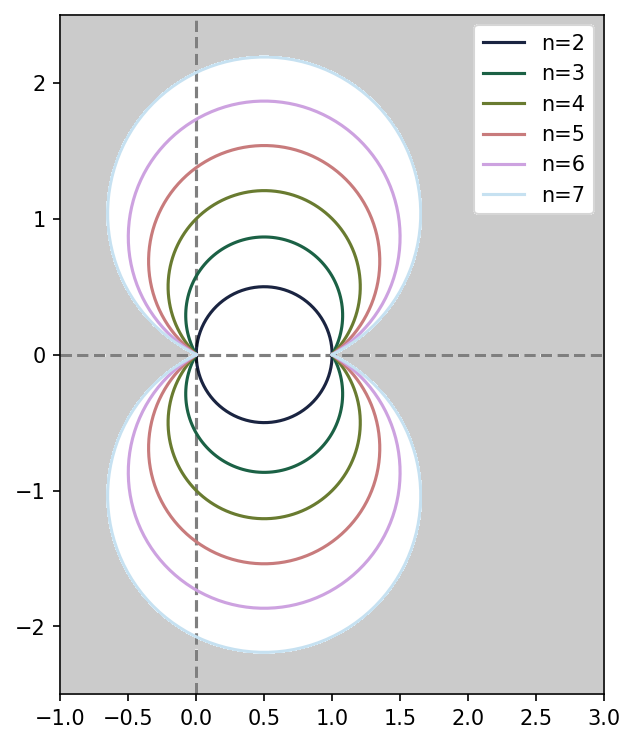}}
\caption{Set of possible eigenvalues of $P_0$ matrices and DPP kernels of size $n\times n$ for $n\in\{2,\cdots, 7\}$.\label{fig: ev}}
\end{figure}

\noindent We illustrate the set of possible values for the eigenvalues of $K$ and $L$ in Figure \ref{fig: ev}.

\noindent Since the distribution of the number of points of a DPP only depends on the eigenvalues of its kernel then a useful method for finding the distribution of $|X|$, when $X\sim\DPP(K)$, is to find a DPP kernel $\tilde K$ with the same eigenvalues as $K$ but with a simpler form from which we can easily infer the distribution of its number of points. Using this method, we first get the following result.

\begin{proposition}\label{prop: Dist_number_point_1}
Let $K$ be a DPP kernel with only real eigenvalues (written with multiplicities) $\lambda_1,\cdots,\lambda_n\in [0, 1]$. If $X\sim\DPP(K)$ then $|X|$ has the same distribution as $B_1+\cdots+B_n$ where $B_1,\cdots, B_n$ are independent Bernoulli random variables with $B_i\sim b(\lambda_i)$.
\end{proposition}
\begin{proof}
The matrix $D(\lambda)$ has the same eigenvalues as $K$ and is a DPP kernel. It corresponds to choosing if each $i\in [n]$ is in $X$ or not independently with probability $\lambda_i$ hence the result.
\end{proof}

\noindent Proposition \ref{prop: Dist_number_point_1} is well-known when $K$ is symmetric but shows that the behavior of $|X|$ as a sum of independent Bernoulli random variables stay consistent across all DPPs whose kernel has only real eigenvalues. Using the same trick, we also recover a result on the behavior of $|X|$ when the kernel eigenvalues can be complex but stay close enough to $\frac{1}{2}$. This result was pointed out in \cite{Brunel18} but proved using a different method.

\begin{proposition}\label{prop: Dist_number_point_2}
Let $K$ be a DPP kernel with eigenvalues in $\mathcal{B}_\mathbb{C}\left(\frac 1 2, \frac 1 2\right)$. We write $\lambda_1,\cdots,\lambda_k$ the real eigenvalues (written with multiplicities) of $K$ and $(\mu_1,\bar\mu_1),\cdots,(\mu_l,\bar\mu_l)$ the complex pairs of eigenvalues (written with multiplicities) of $K$. If $X\sim\DPP(K)$ then $|X|$ has the same distribution as $B_1+\cdots+B_k+C_1+\cdots+C_l$ where $B_i\sim b(\lambda_i)$ and the $C_j$ are random variables in $\{0, 1, 2\}$ satisfying
\[\mathbb{P}(C_j=0)=|\mu_j-1|^2,~~\mathbb{P}(C_j=1)=2\left(\frac 1 4 - \left|\mu_j - \frac 1 2\right|^2\right)~~\mbox{and}~~\mathbb{P}(C_j=2)=|\mu_j|^2.\]
All $B_i$ and $C_j$ are mutually independent from each other.
\end{proposition}
\begin{proof}
We write $\mu_j=a_j+b_ji$ for each $j$ hence $K$ has the same eigenvalues as
\[\tilde K=\begin{pmatrix} \lambda_1 & & 0 & & & 0 \\
& \ddots & & & &\\
0 & & \lambda_k & & &\\
& & & M_1 & & 0\\
& & & & \ddots & \\
0 & & & 0 & & M_l
\end{pmatrix}\]
where each $M_j$ is the $2\times 2$ matrix $\begin{pmatrix} a_j & -b_j \\ b_j & a_j\end{pmatrix}$. It is straightforward to see that each $M_j$ is a DPP kernel by directly computing the associated probabilities. Therefore, by Proposition \ref{prop: Already_proved} \ref{subprop: Block_diag} the matrix $\tilde K$ is a DPP kernel and if $X\sim \DPP(K)$ then $|X|$ has the same distribution has the sum of independent random variables $B_i\sim b(\lambda_i)$ and $C_j$, where $C_j$ is the number of points of a DPP with kernel $M_j$ giving us
\begin{equation*}
\left\{\begin{array}{l}
\mathbb{P}(C_j=0)=\det(I_2-M_j)=(a_j-1)^2+b_j^2=|\mu_j-1|^2;\\
\mathbb{P}(C_j=2)=\det(M_j)=a_j^2+b_j^2=|\mu_j|^2;\\
\mathbb{P}(C_j=1)=1-|\mu_j|^2-|\mu_j-1|^2=1-\frac{|\mu_j-(\mu_j-1)|^2+|\mu_j+(\mu_j-1)|^2}{2}=\frac{1}{2}-\frac{1}{2}|2\mu_j-1|^2.
\end{array}\right.\qedhere
\end{equation*}
\end{proof}

\noindent 
%Proposition \ref{prop: Dist_number_point_1} is a standard result when $K$ is symmetric but we can see that this result stays true for any DPP kernel with real eigenvalues. 
%Proposition \ref{prop: Dist_number_point_2} is a direct extension of \ref{prop: Dist_number_point_1} and can be applied to any DPP kernel of the form $K=\frac 1 2(I_n-M)$ where $\|M\|_2\leq 1$ since, using the properties of the spectral norm, any eigenvalue $\lambda$ of $K$ satisfy $|2\lambda-1|\leq\|M\|_2\leq 1$ hence $\lambda\in \mathcal{B}_\mathbb{C}\left(\frac 1 2, \frac 1 2\right)$. 

Note that Proposition \ref{prop: Dist_number_point_2} can be applied to any DPP kernel of the form $K=\frac 1 2(I_n-M)$ where $\|M\|_2\leq 1$ since, using the properties of the spectral norm, any eigenvalue $\lambda$ of $K$ satisfy $|2\lambda-1|\leq\|M\|_2\leq 1$ hence $\lambda\in \mathcal{B}_\mathbb{C}\left(\frac 1 2, \frac 1 2\right)$. Both Proposition \ref{prop: Dist_number_point_1} and Proposition \ref{prop: Dist_number_point_2} can be further extended using the same trick combined with Proposition \ref{prop: Companion_matrix_K}. Let $\lambda_1,\cdots,\lambda_n$ be the eigenvalues of a DPP kernel $K$ and consider a partition $S_1,\cdots, S_k$ of $[n]$ such that, for each $S_i$, $\prod_{s\in S_i}(X+\frac{\lambda_s}{\lambda_s-1})$ is a polynomial with non-negative coefficients. Then, the distribution of the number of points of a DPP with kernel $K$ is the same as the distribution of $C_1+\cdots+C_k$ where $C_1, \cdots, C_k$ are independent random variables satisfying $C_i\in\{0,\cdots, |S_i|\}$. While such a partition always exist, finding the optimal factorization of a polynomial with non-negative coefficients into a product of lower order polynomials with non-negative coefficients is a difficult problem to our knowledge \cite{Briggs}.

\subsection{Construction with some specific eigenvalues}

\noindent It was shown in \cite{Kellogg} how to construct some specific $P_0$ matrices with a given set of eigenvalues. By adapting these results into the DPP setting we can also construct some DPP kernels with a given set of eigenvalues.

\begin{proposition}
Let $\lambda_1,\cdots,\lambda_n\in\mathbb{C}$ be the eigenvalues (taken with multiplicities) of some $P_0$ matrix $\tilde L$. We consider the polynomial 
\[P(X)=\prod_{i=1}^n (X+\lambda_i)=c_0+c_1X+\cdots+c_{n-1}X^{n-1}+X^n\]
and the kernel
\begin{equation}\label{eq: Companion_matrix}
L = \begin{pmatrix} 0 & -1 & 0 & \cdots & 0\\
\vdots & \ddots & \ddots & \ddots & \vdots\\
\vdots &  & \ddots & \ddots & 0 \\
0 & \cdots & \cdots & 0 & -1 \\
c_0 & \cdots & \cdots & c_{n-2} & c_{n-1}
\end{pmatrix}.
\end{equation}
Then $L$ is also a $P_0$ matrix. Moreover, if $X\sim \DPP_L(L)$ then
\[\mathbb{P}(X=S)=\left\{\begin{array}{l}
\displaystyle\frac{c_{k-1}}{c_0+\cdots+c_{n-1}+1}~~\mbox{if}~S=\{k,\cdots, n\};\\
\displaystyle\frac{1}{c_0+\cdots+c_{n-1}+1}~~\mbox{if}~S=\emptyset;\\
0~~\mbox{otherwise}.
\end{array}
\right.\]
\end{proposition}
\begin{proof}
Using standard results on characteristic polynomials we have $c_k=e_k(\lambda_1,\cdots, \lambda_n)$ which is equal to the sum of all principal minors of $\tilde L$ of size $k$ \cite[Theorem 1.2.16]{Horn} and is therefore non-negative since $\tilde L$ is a $P_0$ matrix. Moreover, $-L$ is the companion matrix of $P(X)$ hence the eigenvalues of $L$ are the root of $P(-X)$ corresponding to the $\lambda_i$. Now, let $S\subset [n]$. Obviously, if $n\notin S$ then $L_S$ is a submatrix of a triangular matrix whose diagonal entries are all zeros hence $\det(L_S)=0$. Otherwise, we write $S=\{s_1,\cdots,s_{k-1},s_k\}$ with $s_1<\cdots<s_k=n$ and we can write $L_S$ as
\[L_S = \begin{pmatrix} 0 & & & & \\
\vdots & & M & \\
0 & & & \\
c_{s_1} & c_{s_2} & \cdots & c_{s_k}
\end{pmatrix}.\]
where $M$ is an upper triangular matrix whose diagonal elements are $L_{s_1,s_2},\cdots,L_{s_{k-1}, s_k}$ and therefore 
\[\det(L_S)=(-1)^{k+1}c_{s_1}\det(M)=(-1)^{k+1}c_{s_1}\prod_{i=1}^{k-1}L_{s_{i-1}, s_i}.\]
Now, $L_{s_{i-1}, s_i}$ is equal to $-1$ if $s_i=s_{i-1}+1$ and $0$ otherwise hence $\det(L_S)$ is nonzero only when $S$ is of the form $\{k,\cdots, n\}$ with $\det(L_S)=c_k\geq 0$ in this case. We can conclude that $L$ is a $P_0$ matrix and we get the expression of $\bP(X=S)$ from \eqref{eq: L_prob}.
\end{proof}

\noindent As a direct consequence, computing $K=I_n-(I_n+L)^{-1}$ when $L$ is of the form \eqref{eq: Companion_matrix} yields the following result.

\begin{proposition}\label{prop: Companion_matrix_K}
Let $\lambda_1,\cdots,\lambda_n\in\mathbb{C}\backslash\{1\}$ be the eigenvalues of some DPP kernel $\tilde K$. We consider the polynomial 
\[P(X)=\prod_{i=1}^n \left(X+\frac{\lambda_i}{1-\lambda_i}\right)=c_0+c_1X+\cdots+c_{n-1}X^{n-1}+X^n\]
and the kernel
\[K = \frac{1}{c_0+\cdots+c_{n-1}+1}\begin{pmatrix} c_0 & c_0+c_1 & \cdots & c_0+\cdots+c_{n-1}\\
\vdots & \vdots & & \vdots\\
c_0 & c_0+c_1 & \cdots & c_0+\cdots+c_{n-1}\\
\end{pmatrix} - \begin{pmatrix} 0  & 1 & \cdots & 1\\
\vdots & \ddots & \ddots & \vdots \\
\vdots & & \ddots & 1\\
0 & \cdots & \cdots & 0
\end{pmatrix}.\]
Then $K$ is also a DPP kernel. Moreover, if $X\sim \DPP(K)$ then
\[\mathbb{P}(X=S)=\left\{\begin{array}{l}
\displaystyle\frac{c_{k-1}}{c_0+\cdots+c_{n-1}+1}~~\mbox{if}~S=\{k,\cdots, n\};\\
\displaystyle\frac{1}{c_0+\cdots+c_{n-1}+1}~~\mbox{if}~S=\emptyset;\\
0~~\mbox{otherwise}.
\end{array}
\right.\]
\end{proposition}

\subsection{Examples of DPPs with nonsymmetric kernels}

The literature already contains quite a few examples of DPPs with generic kernels. Numerous examples of determinantal point processes are given in \cite[Section 2.5]{Soshnikov} with some of them having nonsymmetric kernels. We mention for example the Coulomb gas in \cite{Cornu} or the result in \cite{Borodin99} showing that, when taking a random partition using a Plancherel measure (later generalized in \cite{Okounkov} to Schur measures) and transforming it into a point process on $\frac{1}{2}+\mathbb{Z}$ using a natural transformation we get a DPP with a nonsymmetric kernel. DPP kernels $K$ satisfying $K_{i, j}= \pm K_{j, i}$ have been studied in \cite{Brunel18} and L-ensemble kernels of the form $L=VV^T + B(D-D^T)B^T$ are considered in \cite{Gartrell21, Han, Han22}. We also mention the result of \cite{Borodin09} that $1$-dependent binary processes are DPPs whose kernel can be written as an upper Hessenberg matrix with various examples given in the aforementioned paper. We complete this list by giving a few additional simple examples of DPPs with nonsymetric kernels.

\subsubsection*{Rank one matrix}

For rank one matrices it is easy to characterize which one are DPP kernels or not.

\begin{proposition}
Let $K\in\mathcal{M}_n(\bR)$ be of rank $1$. We write $K=\lambda uv^T$ where $\lambda$ is the only non-zero eigenvalue of $K$ and $u,v\in\bR^n$ are left and right eigenvectors chosen such that $\langle u, v\rangle=1$. Then $K$ is a DPP kernel if and only if 
\[\forall i\in [n],~u_iv_i\geq 0~\mbox{and}~\lambda\in [0, 1].\]
\end{proposition}
\begin{proof}
Since $K$ is of rank one then $\det(K_S)=0$ when $|S|\geq 2$. Denoting by $\mu$ the determinantal measure associated with $K$ we then get $\mu(\{S\})=0$ if $|S|\geq 2$. If $S=\{i\}$ for some $i\in [n]$ then $\mu(\{i\})=\sum_{S\ni i} \mu(S)= K_{i, i} = \lambda u_iv_i$ and thus $\mu(\emptyset)= 1 - \lambda \sum_{i=1}^n u_iv_i = 1 - \lambda$, proving the result.
\end{proof}

\subsubsection*{Rank one perturbation of the half identity matrix}

Similarly, we can also characterize exactly which rank one perturbations of $\frac{1}{2}I_n$ are DPP kernels.

\begin{proposition}
Let $K\in\mathcal{M}_n(\bR)$ be of the form $K=\frac{1}{2}(I_n - uv^T)$ for some vectors $u,v\in\bR^n$. Then $K$ is a DPP kernel if and only if $\sum_{i=1}^n|u_iv_i|\leq 1$.
\end{proposition}
\begin{proof}
Let $S\subset [n]$. Using identity \eqref{eq: Sum_minors} we have
\begin{multline*}
\det(D(\cara{S})(I_n-K)+D(\cara{S^c})K) = \frac{1}{2^n}\det(I_n + D(\cara{S} - \cara{S^c})uv^T) \\
= \frac{1}{2^n}\sum_{T\subset [n]}\det((D(\cara{S} - \cara{S^c})uv^T)_T)
= \frac{1}{2^n}\left(1 + \sum_{i\in S} u_iv_i - \sum_{i\notin S} u_iv_i\right),
\end{multline*}
where the last identity is a consequence of $D(\cara{S} - \cara{S^c})uv^T$ being a rank one matrix and thus having vanishing principal minors of size $\geq 2$. Thus, the lowest possible value of $\det(D(\cara{S})(I_n-K)+D(\cara{S^c})K)$ among all $S\subset [n]$ is $\frac{1}{2^n}(1 - \sum_{i=1}^n|u_iv_i|)$ and we conclude using Theorem \ref{theo: Caracterisation_DPPs} \ref{cara1}.
\end{proof}

\noindent Note that, as a consequence of Corollary \ref{corr: Alternate_LL}, we get the identity
\[\mathbb{P}(X=S)=\frac{1}{2^n}\left(1-\sum_{i\in S}u_iv_i+\sum_{i\not\in S}u_iv_i\right)~\Rightarrow~\mathbb{P}(X=S)+\mathbb{P}(X=S^c)=\frac{1}{2^{n-1}}\]
giving the following characterization:

\begin{corollary}\label{prop: Simu_half_identity_rank_1}
Let $u,v\in\bR^n$ such that $\sum_{i=1}^n|u_iv_i|\leq 1$. Let $S$ be a random subset of $[n]$ with a uniform distribution on $\mathcal{P}([n])$. We define the random subset $X$ of $[n]$ conditionally to $S$ by
\[\forall T\subset [n],~\mathbb{P}(X=T|S)=\left\{\begin{array}{l}
\frac{1}{2}\left(1-\sum_{i\in S}u_iv_i+\sum_{i\not\in S}u_iv_i\right)~\mbox{if}~T=S;\\
\frac{1}{2}\left(1-\sum_{i\not\in S}u_iv_i+\sum_{i\in S}u_iv_i\right)~\mbox{if}~T=S^c;\\
0~\mbox{otherwise.}
\end{array}
\right.\]
Then $X$ is a DPP with kernel $K=\frac{1}{2}(I_n - uv^T)$.
\end{corollary}

\subsubsection*{Uniform random set with a given parity}

The symmetric DPP kernel $\frac{1}{2}I_n$ corresponds to the uniform distribution on the $2^n$ subsets of $\cP([n])$. Interestingly, when considering instead the uniform distribution on the subsets of $\cP([n])$ whose cardinal is even (or odd) then we still get a DPP but with a nonsymmetric kernel.

\begin{proposition}
Let $X$ be uniformly distributed on $\cP([n])$. Then
\begin{equation}\label{eq: unif_cond_even}
X{\Big |} |X|~\mbox{even}\sim\DPP\left(\frac 1 2\begin{pmatrix} 1 & 0 & \cdots & 0 & -1\\
1 & \ddots & \ddots & & 0\\
0 & \ddots & \ddots & \ddots & \vdots\\
\vdots & \ddots & \ddots & \ddots & 0\\
0 & \cdots & 0 & 1 & 1
\end{pmatrix}\right)
\end{equation}
\begin{equation}\label{eq: unif_cond_odd}
X{\Big |} |X|~\mbox{odd}\sim\DPP\left(\frac 1 2\begin{pmatrix} 1 & 0 & \cdots & 0 & 1\\
1 & \ddots & \ddots & & 0\\
0 & \ddots & \ddots & \ddots & \vdots\\
\vdots & \ddots & \ddots & \ddots & 0\\
0 & \cdots & 0 & 1 & 1
\end{pmatrix}\right)
\end{equation}
\end{proposition} 
\begin{proof}
The kernel in both identity \eqref{eq: unif_cond_even} and identity \eqref{eq: unif_cond_odd} is of the form $\frac{1}{2}(I_n - M)$ where $M$ is an orthogonal matrix satisfying $M_{i, j} \in \{-1, 1\}$ when $i - j \equiv 1 [n]$ and $0$ otherwise. These are indeed DPP kernels by Proposition \ref{prop: perturb_half_identity} and they are even a special case of what we, later on, refer to as an orthogonal DPP (see Section \ref{sec: Mixing}). Moreover, we have $\det(M_S)=0$ when $S\neq [n]$ in both cases, $\det(M)=(-1)^n$ for \eqref{eq: unif_cond_even} and $\det(M)=(-1)^{n + 1}$ for \eqref{eq: unif_cond_odd}. Letting $X\sim\DPP(K)$ where $K$ is the kernel in either \eqref{eq: unif_cond_even} or \eqref{eq: unif_cond_odd} and using \eqref{eq: Sum_minors} then give, for all $S\neq [n]$,
\[\bP(S\subset X) = \det\left(\frac{1}{2}(I_{|S|} - M_S)\right)=\frac{1}{2^{|S|}}\sum_{T\subset S}\det(M_T) = \frac{1}{2^{|S|}}\]
and
\[\bP([n]\subset X) = \left\{\begin{array}{rl} \frac{1}{2^n}\big(1 + (-1)^n\big)~&\mbox{if}~K~\mbox{is from}~\eqref{eq: unif_cond_even},\\
\frac{1}{2^n}\big(1 + (-1)^{n + 1}\big)~&\mbox{if}~K~\mbox{is from}~\eqref{eq: unif_cond_odd}.
\end{array}\right.\]
This corresponds to the inclusion probabilities of a uniform random set conditioned to have an even (resp. odd) number of points for the kernel in identity \eqref{eq: unif_cond_even} (resp. \eqref{eq: unif_cond_odd}), concluding the proof.
\end{proof}

\noindent We note that this result can also be seen as a special case of a $1$-dependent binary process as in \cite{Borodin09}.

\section{Properties of general DPPs}\label{sec: prop_DPPs}

\noindent We saw in Proposition \ref{prop: star-shaped} that $\frac{1}{2}I_n$ can be seen as the center of the set of DPP kernels and we showed in Proposition \ref{prop: perturb_half_identity} that the set of DPP kernels contains the ball centered at $\frac{1}{2}I_n$ with radius $\frac{1}{2}$ for the spectral norm. This suggests that a natural way of writing DPP kernels is as $K=\frac{1}{2}(I_n-M)$, giving rise to a new kernel $M$. We focus in this section on the properties of this kernel $M$ and show that we can use it to extend some of the classical properties of DPPs with symmetric kernels.

\subsection{The parity kernel}

We begin by pointing out that the relationship between $M$ and the likelihood kernel $L$ is a standard matrix transformation called the Cayley transform whose main property is that it maps skew-symmetric matrices to special orthogonal matrices.

\begin{proposition}\label{prop: Cayley}
Let $K=\frac{1}{2}(I_n-M)$ be a DPP kernel such that $I_n+M$ is invertible. Then $L$ is the Cayley transform of $M$:
\[L=(I_n-M)(I_n+M)^{-1}.\]
\end{proposition}

\noindent We recall that the principal minors of the kernels $K$ and $L$ gives, respectively, the inclusion and exact probabilities of the DPP. This is why $K$ and $L$ are often called, respectively, the correlation and likelihood kernel of a DPP. In comparison, the principal minors of the kernel $M$ also have a nice direct interpretation as a consequence of the Laplace transform formula (see Proposition \ref{prop: Already_proved}~\ref{subprop: Laplace}):
\begin{equation}\label{eq: defin_laplace_M}
X\sim \DPP(K)~\Rightarrow~\forall S\subset [n],~\bE\left[(-1)^{|X\cap S|}\right]=\det((I_n-2K)_S)=\det(M_S)
\end{equation}
hence
\[\mathbb{P}(|X\cap S|~\mbox{is even})= \frac{1}{2}(1+\det(M_S))~\mbox{and}~\mathbb{P}(|X\cap S|~\mbox{is odd})= \frac{1}{2}(1-\det(M_S)).\]
The principal minors of the kernel $M$ thus convey information about the parity of the number of points of a DPP in a given subset. For this reason we decide to call $M$ the parity kernel of the DPP. 
\begin{definition}
Let $X$ be a DPP with kernel $K\in\cM_n(\bR)$. We define the matrix $M=I_n - 2K$ as the \textbf{parity kernel} of the DPP.
\end{definition}
\noindent Generalizing \eqref{eq: defin_laplace_M} gives the following proposition.
\begin{proposition}\label{prop: General_parity}
Let $S_1,\cdots, S_k$ be disjoint subsets of $[n]$. Let $R_1,\cdots, R_k$ defined by
\[R_i=(-1)^{|X\cap S_i|}=\left\{\begin{array}{l} 1~\mbox{if}~|X\cap S_i|~\mbox{is even,}\\
-1~\mbox{if}~|X\cap S_i|~\mbox{is odd.}\end{array}
\right.\]
Then, for all $\epsilon\in\{-1, +1\}^k$,
\[\mathbb{P}(\forall i,~R_i=\epsilon_i)=\frac{1}{2^k}\sum_{T\subset [k]}\left(\prod_{i\in T}\epsilon_i\right)\det(M_{\cup_{i\in T}S_i})\]
\end{proposition}
\begin{proof}
Using expression \eqref{eq: defin_laplace_M} and the fact that $1+\epsilon_iR_i$ is equal to $0$ if $R_i\neq \epsilon_i$ and $2$ if $R_i = \epsilon_i$ we get
\begin{align*}
\mathbb{P}(\forall i,~R_i=\epsilon_i)&=\frac{1}{2^k}\bE\left[\prod_{i=1}^k\left(1+\epsilon_iR_i\right)\right]\\
&=\frac{1}{2^k}\bE\left[\prod_{i=1}^k\left(1+\epsilon_i(-1)^{|X\cap S_i|}\right)\right]\\
&=\frac{1}{2^k}\sum_{T\subset [k]}\bE\left[\prod_{i\in T}\epsilon_i(-1)^{|X\cap S_i|}\right]\\
&=\frac{1}{2^k}\sum_{T\subset [k]}\left(\prod_{i\in T}\epsilon_i\right)\bE\left[(-1)^{|X\cap (\cup_{i\in T}S_i)|}\right]\\
&=\frac{1}{2^k}\sum_{T\subset [k]}\left(\prod_{i\in T}\epsilon_i\right)\det(M_{\cup_{i\in T}S_i}).
\end{align*}
\end{proof}

\noindent As a consequence, for any $S\subset [n]$ we can write
\[\left\{\begin{array}{l}
\bP(|X\cap S|~\mbox{even and}~|X\cap S^c|~\mbox{even})=\frac{1}{4}\left(1+\det(M_S)+\det(M_{S^c})+\det(M)\right);\\
\bP(|X\cap S|~\mbox{even and}~|X\cap S^c|~\mbox{odd})=\frac{1}{4}\left(1+\det(M_S)-\det(M_{S^c})-\det(M)\right);\\
\bP(|X\cap S|~\mbox{odd and}~|X\cap S^c|~\mbox{even})=\frac{1}{4}\left(1-\det(M_S)+\det(M_{S^c})-\det(M)\right);\\
\bP(|X\cap S|~\mbox{odd and}~|X\cap S^c|~\mbox{odd})=\frac{1}{4}\left(1-\det(M_S)-\det(M_{S^c})+\det(M)\right).
\end{array}\right.\]

\noindent Additionally, there is an almost symmetric relationship between the probabilities of a DPP and the minors of $M$ due to the Cayley involution (Proposition \ref{prop: Cayley}) mapping the kernel $M$ to the kernel $L$.

\begin{proposition}
Let $K=\frac{1}{2}(I_n-M)\in\mathcal{M}_n(\bR)$ be a DPP kernel and $X\sim \DPP(K)$. Then, for any $S\subset [n]$:
\begin{equation}\label{eq: Proba_to_M}
\mathbb{P}(X=S)=\frac{1}{2^n}\sum_{T \subset [n]}(-1)^{|T \cap S|}\det(M_T);
\end{equation}
\begin{equation}\label{eq: M_to_Proba}
\det(M_S)=\sum_{T \subset [n]}(-1)^{|T \cap S|}\mathbb{P}(X=T).
\end{equation}
\end{proposition}
\begin{proof}
Identity \eqref{eq: M_to_Proba} is an immediate consequence of \eqref{eq: defin_laplace_M} since
\[\bE\left[(-1)^{|X\cap S|}\right]=\sum_{T \subset [n]}(-1)^{|T \cap S|}\mathbb{P}(X=T),\]
while \eqref{eq: Proba_to_M} is an immediate consequence of Proposition \ref{prop: General_parity} since $X=S$ is equivalent to $|X\cap\{i\}|$ being even for $i\notin S$ and odd for $i\in S$.
\end{proof}

\subsection{Mixing of orthogonal DPPs}\label{sec: Mixing}
The most common method used to simulate DPPs with a symmetric kernel $K$ is based on a first result stating that all the eigenvalues $\lambda_i$ of $K$ can be replaced with independent Bernoulli random variables with parameter $\lambda_i$ (see \cite[Theorem 4.5.3]{Hough} for example), turning $K$ into a projection matrix. Such a DPP is called a \textbf{projection DPP} and its various nice properties (such as having a deterministic number of points) allows for a simple simulation method, often called the HKPV algorithm \cite[Algorithm 4.4.2]{Hough}, that can be used for both finite and continuous DPPs. We generalize this principle for kernels of the form $K=\frac{1}{2}(I_n - M)$ with $\|M\|_2\leq 1$ by proving that we can replace the singular values of the parity kernel $M$ by random variables with values in $\{-1, +1\}$. The kernel $M$ then become an orthogonal matrix and we call the associated determinantal distribution an orthogonal DPP.

\begin{definition}
We call a DPP $X$ whose parity kernel $M$ is an orthogonal matrix an \textbf{orthogonal DPP}.
\end{definition}

\begin{proposition} \label{prop: bad_mixing}
Let $M\in\mathcal{M}_n(\bR)$ be a matrix such that $\|M\|_2\leq 1$ and write $M=PD(\sigma)Q^T$ its singular value decomposition where $P,Q$ are $n\times n$ orthogonal matrices and $\sigma=(\sigma_1,\cdots,\sigma_n)\in[0,1]^n$ is the vector of singular values of $M$. Let $B=(B_1,\cdots,B_n)$ be a vector of $n$ independent random variables taking value in $\{-1, 1\}$ and satisfying $\bE[B_i]=\sigma_i$. Define $\tilde K = \frac{1}{2}(I_n - PD(B)Q^T)$. Note that $\tilde K$ is always a DPP kernel by Proposition \ref{prop: perturb_half_identity}. If $X$ is a DPP with kernel $\tilde K$ conditionally to $B$ then $X$ is a DPP with kernel $K=\frac{1}{2}(I_n - M)$.
\end{proposition}
\begin{proof}
Let $S\subset [n]$. Using \eqref{eq: Sum_minors} we can write
$$\det(K_S)=\frac{1}{2^{|S|}}\det(I_{|S|} - M_S)=\frac{1}{2^{|S|}}\sum_{T\subset S}\det(-M_T)=\frac{1}{2^{|S|}}\sum_{T\subset S}(-1)^{|T|}\det(P_{T, [n]}D(\sigma)Q_{[n], T}^T).$$
Using the Cauchy-Binet formula twice and the fact that $\det(D(\sigma)_{U,V})=(\prod_{i\in U}\sigma_i)\cara{U=V}$ gives 
\begin{align*}
\det(K_S)&=\frac{1}{2^{|S|}}\sum_{T\subset S}(-1)^{|T|}\sum_{\substack{U,V\subset [n]\\ |U|=|V|=|T|}}\det(P_{T,U})\det(D(\sigma)_{U,V})\det(Q_{T, V})\\
&=\frac{1}{2^{|S|}}\sum_{T\subset S}(-1)^{|T|}\sum_{\substack{U\subset [n]\\ |U|=|T|}}\det(P_{T,U})\det(Q_{T, U})\prod_{i\in U}\sigma_i.
\end{align*}
With the same reasoning we also get
\[\det(\tilde K_S)=\frac{1}{2^{|S|}}\sum_{T\subset S}(-1)^{|T|}\sum_{\substack{U\subset [n]\\ |U|=|T|}}\det(P_{T, U})\det(Q_{T, U})\prod_{i\in U}B_i,\]
hence $\bE[\det(\tilde K_S)]=\det(K_S)$ and thus
\begin{equation*}
\mathbb{P}(S\subset X)=\bE[\mathbb{P}(S\subset X| B_1,\cdots, B_n)]=\bE[\det(\tilde K_S)]=\det(K_S).\qedhere
\end{equation*}
\end{proof}

\subsection{Properties of orthogonal DPPs}

We point out that a slight potential source of confusion is that the name "projection DPP" refers to the correlation kernel $K$ being a projection matrix, while the name "orthogonal DPP" refers to the parity matrix $M$ being orthogonal. The kernel $K$ of orthogonal DPPs satisfy
\begin{equation}\label{eq: K_orthoDPP}
K^TK = KK^T = \frac{1}{2}(K + K^T),
\end{equation}
which does not corresponds to a common class of matrices to our knowledge, except that it is a particular case of normal matrices. When $K$ is symmetric, identity \eqref{eq: K_orthoDPP} becomes $K^2=K$ meaning that orthogonal DPPs generalize projection DPPs. Just like projection DPPs we show that orthogonal DPPs are stable by conditioning with respect to the presence of a point.

\begin{proposition}
Let $X$ be an orthogonal DPP with associated kernel $K=\frac{1}{2}(I_n - M)$, $M$ being and orthogonal matrix. Let $i\in [n]$ such that $K_{i,i}\neq 0$. Then, the distribution of $X\backslash \{i\}$ conditionally to $i\in X$ is an orthogonal DPP.
\end{proposition}
\begin{proof}
Let $\tilde K$, defined as in \eqref{eq: conditional_K}, be the kernel of the distribution of $X\backslash \{i\}$ conditionally to $i\in X$ and define $\tilde M\defeq I_n - 2\tilde K$. Since $K$ satisfy \eqref{eq: K_orthoDPP} then
\[K(K_{i, [n]})^T=\frac{1}{2}(K_{[n], i}+(K_{i, [n]})^T)~\Rightarrow~M(K_{i, [n]})^T=(K_{i, [n]})^T-(K_{[n], i}+(K_{i, [n]})^T)=-K_{[n], i}\]
and
\[\langle (K_{i, [n]})^T, (K_{i, [n]})^T\rangle=\frac{1}{2}(K_{i,i}+K_{i,i})=K_{i, i}.\]
Hence, if we define the unit vector $v=\frac{1}{\|(K_{i, [n]})^T\|}(K_{i, [n]})^T$ then
\[\tilde M = I_n - 2K + \frac{2}{K_{i,i}}K_{[n], i}K_{i, [n]}=M - \frac{2}{\|(K_{i, [n]})^T\|^2}M(K_{i, [n]})^TK_{i, [n]}=M(1-2vv^T).\]
Therefore, $\tilde M$ is the product of an orthogonal matrix and an householder reflection and thus an orthogonal matrix itself.
\end{proof}

\noindent More generally we show that, unlike projection DPPs, orthogonal DPPs are also stable by particle-hole involution.

\begin{proposition}\label{prop: ortho_stable_particle_hole}
Let $X$ be an orthogonal DPP with associated kernel $K=\frac{1}{2}(I_n - M)$, $M$ being and orthogonal matrix. Let $S\subset [n]$ and $\tilde X= (X\cap S^c)\cup(X^c\cap S)$. Then, $\tilde X$ is an orthogonal DPP.
\end{proposition}
\begin{proof}
We saw in Proposition \ref{prop: Particle_Hole} that $\tilde X$ is a DPP with kernel $\tilde K$ satisfying identity \eqref{eq: PH_K}. Thus,
\begin{align*}
\tilde M &= I_n - 2\tilde K\\
&= I_n - 2D(\cara{S})(I_n-K) -2(I_n - D(\cara{S}))K\\
&= I_n - D(\cara{S})(I_n + M) -(I_n - D(\cara{S}))(I_n - M)\\
&= I_n - D(\cara{S}) -D(\cara{S})M -I_n + M + D(\cara{S}) - D(\cara{S})M\\
&= (I_n - 2D(\cara{S}))M
\end{align*}
Since $D(\cara{S})$ is a symmetric projection matrix then $(I_n-2D(\cara{S}))M$ is orthogonal.
\end{proof}

\noindent Note that a direct consequence of the proof of Proposition \ref{prop: ortho_stable_particle_hole} is that the particle involution with respect to a singleton $\{i\}$ multiply $M$ by the householder reflection with respect to the $i$-th canonical vector.

\medskip

\noindent Finally, an immediate consequence of the decomposition of real orthogonal matrices gives the following decomposition of the $K$ and $L$ kernels of orthogonal DPPs.

\begin{proposition}\label{prop: decomp_K_orthogonal}
Let $K$ be the kernel of an orthogonal DPP. There exists an $n\times n$ orthogonal matrix $P$ such that $PKP^T$ is written as
\begin{equation}\label{eq: decomp_K_orthogonal}
\begin{pmatrix}
R_1 & & 0 & & & \\
 & \ddots & & & 0 & \\
 0 & & R_k & & & \\
 & & & \delta_1 & & 0 \\
 & 0 & & & \ddots & \\
 & & & 0 & & \delta_{n-2k}
\end{pmatrix}
\end{equation}
where $\delta_i\in\{0, 1\}$ and each $R_i$ is a $2\times 2$ matrix of the form
\[R_i=\begin{pmatrix}
\frac{1-\cos(\theta_i)}{2} & -\frac{\sin(\theta_i)}{2}\\
\frac{\sin(\theta_i)}{2} & \frac{1-\cos(\theta_i)}{2}
\end{pmatrix},~\theta_i\in (0,\pi).\]
If all $\delta_i$ are equal to $0$ then $L$ is a well defined skew-symmetric matrix and $PLP^T$ is written
\begin{equation}\label{eq: decomp_L_orthogonal}
\begin{pmatrix}
R'_1 & & 0 & & & \\
 & \ddots & & & 0 & \\
 0 & & R'_k & & & \\
 & & & 0 & & 0 \\
 & 0 & & & \ddots & \\
 & & & 0 & & 0
\end{pmatrix}~~\mbox{where}~~R'_i=\begin{pmatrix}
0 & \tan\left(\frac{\theta_i}{2}\right)\\
-\tan\left(\frac{\theta_i}{2}\right) & 0
\end{pmatrix}.
\end{equation}
\end{proposition}

\noindent We point out that a DPP with kernel $K=\begin{pmatrix}\frac{1-\cos(\theta)}{2} & -\frac{\sin(\theta)}{2}\\
\frac{\sin(\theta)}{2} & \frac{1-\cos(\theta)}{2}
\end{pmatrix}$ satisfy
\[\mathbb{P}(X=\emptyset)=\frac{1+\cos(\theta)}{2},~\mathbb{P}(X=\{1\})=\mathbb{P}(X=\{2\})=0~\mbox{and}~\mathbb{P}(X=\{1, 2\})=\frac{1-\cos(\theta)}{2},\]
telling us that the cardinal of orthogonal DPPs with the decomposition \eqref{eq: decomp_K_orthogonal} is of the form $l+2B$ where $l$ is the number of $\delta_i$ equal to $1$ and $B$ is a Poisson-Binomial distribution with probabilities $\frac{1-\cos(\theta_1)}{2},\cdots,\frac{1-\cos(\theta_k)}{2}$. If $k=0$ then $K$ is a projection matrix and we recover the usual result that projection DPPs have a deterministic amount of points.

\subsection{Simulation of orthogonal DPPs}\label{subsec: DPP_sim}

By Proposition \ref{prop: decomp_K_orthogonal}, the kernel $K$ of an orthogonal DPP has a number $k$ of pairs of complex eigenvalues, a number $l$ of eigenvalues equal to $1$ and a number $n-l-2k$ of eigenvalues equal to $0$. When $k=0$ this corresponds to a projection DPP for whom we already know an efficient sampling algorithm \cite[Algorithm 4.4.2]{Hough}. When $k=1$ we can show that the corresponding orthogonal DPP is a mixture of two projection DPPs.

\begin{proposition}
Let $P$ be a matrix with size $n \times (l+2)$ satisfying $P^TP=I_{l+2}$ and $K$ be the kernel of an orthogonal DPP such that
\begin{equation}\label{eq: decomp_K_orthogonal_k=1}
K=PDP^T~\mbox{where}~D=\begin{pmatrix}
\alpha & -\beta & 0 & \cdots & 0\\
\beta & \alpha & 0 & \cdots & 0\\
0 & 0 & 1 & & 0 \\
\vdots & \vdots & & \ddots & \\
0 & 0 & 0 & & 1
\end{pmatrix}
\end{equation}
with $\alpha\in ]0, 1[$ and $\alpha^2+\beta^2=\alpha$. Now, define $X$ as a projection DPP with kernel $K^{(1)}=PP^T$ and $Y$ as a projection DPP with kernel $K^{(2)}=P_{[n],\{3,\cdots, l+2\}}P_{[n],\{3,\cdots, l+2\}}^T$. Let $Z$ be a random variable defined conditionally to $(X, Y)$ such that $Z=X$ with probability $\alpha$ and $Z=Y$ with probability $1-\alpha$ (the choice being independent from $X$ and $Y$). Then, $Z$ is a DPP with kernel $K$.
\end{proposition}
\begin{proof}
Let $S\subset [n]$. If $|S|\geq l+2$ then $\det(K_S)=0$ since $K$ is of rank $l+2$. Otherwise, the Cauchy-Binet formula applied twice gives
\[\det(K_S)=\sum_{\substack{T, U\subset [l+2]\\ |T|=|U|=|S|}}\det(P_{S, T})\det(D_{T, U})\det(P_{S, U}).\]
Due to the shape \eqref{eq: decomp_K_orthogonal_k=1} of the matrix $D$, $\det(D_{T, U})$ vanishes when $T\cap\{3,\cdots,l+2\}\neq U\cap\{3,\cdots,l+2\}$ or $|T\cap\{1, 2\}|\neq|U\cap\{1, 2\}|$ thus
\begin{multline*}
\det(K_S)=\sum_{\substack{V\subset \{3,\cdots,l+2\}\\ |V|=|S|}}\det(P_{S, V})^2 + \sum_{\substack{V\subset \{3,\cdots,l+2\}\\ |V|=|S|-2}}(\alpha^2+\beta^2)\det(P_{S, V\cup\{1, 2\}})^2\\
+ \sum_{\substack{V\subset \{3,\cdots,l+2\}\\ |V|=|S|-1}}\left(\alpha\det(P_{S, V\cup\{1\}})^2 + \alpha\det(P_{S, V\cup\{2\}})^2+ (\beta-\beta)\det(P_{S, V\cup\{2\}})\det(P_{S, V\cup\{1\}})\right).
\end{multline*}
Using the identity $\alpha^2+\beta^2=\alpha$ and once again the Cauchy-Binet formula we finally get
\begin{align*}
\det(K_S)&=(1-\alpha)\sum_{\substack{V\subset \{3,\cdots,l+2\}\\ |V|=|S|}}\det(P_{S, V})^2 + \alpha\sum_{\substack{V\subset [l+2]\\ |V|=|S|}}\det(P_{S, V})^2\\
&=(1-\alpha)\det(K^{(2)}_{S}) + \alpha\det(K^{(1)}_{S})\\
&=(1-\alpha)\mathbb{P}(S\subset Y) + \alpha\mathbb{P}(S\subset X) = \mathbb{P}(S\subset Z),
\end{align*}
concluding the proof.
\end{proof}

\noindent Unfortunately, this result does not generalize well when $k\geq 2$. To illustrate the issue, let $k\geq 2$ and $l=0$ so that the kernel $L$ is well-defined and of the form \eqref{eq: decomp_L_orthogonal}. We write the $2\times 2$ diagonal blocks in \eqref{eq: decomp_L_orthogonal} as $R'_i=\begin{pmatrix} 0 & \mu_i \\ -\mu_i & 0\end{pmatrix}$
to simplify the notations and we also write $p_i=\frac{\mu_i^2}{1+\mu_i^2}\in]0, 1[$ and $q_i=1-p_i$. For odd-sized sets $S\subset [n]$ we have $\det(L_S)=0$ since $L$ is skew-symmetric. When $S$ is of even size then the Pfaffian minor-summation formula \cite{Ishikawa} gives
\[\det(L_S)=\Pf(L_S)^2=\left(\sum_{\substack{T\subset [n]\\ |T|=|S|}}\Pf(D_T)\det(P_{S, T})\right)^2\]
For $\Pf(D_T)$ to be non-zero then we obviously need to have $|T|\leq 2k$ and $T$ has to satisfy either $T\cap\{i, i+1\}=\emptyset$ or $\{i, i+1\}\subset T$ for any odd $i\leq 2k-1$.
Now, for any set $T\subset [k]$ we write $\phi(T)\defeq\cup_{i\in T}\{2i-1, 2i\}$. Then,
\begin{align*}
\frac{\det(L_S)}{\det(I_n+L)}&=\prod_{i=1}^k(1+\mu_i^2)^{-1}\left(\sum_{\substack{T\subset [k]\\ |T|=|S|/2}}\prod_{i\in T}\mu_i\det(P_{S, \phi(T)})\right)^2\\
&=\left(\sum_{\substack{T\subset [k]\\ |T|=|S|/2}}\prod_{i\in T}\sqrt{p_i}\prod_{i\notin T}\sqrt{q_i}\det(P_{S, \phi(T)})\right)^2.
\end{align*}
When $k=1$ we recover the result of Proposition \ref{eq: decomp_K_orthogonal_k=1}:
\[\frac{\det(L_S)}{\det(I_n+L)}=\left\{\begin{array}{l}
p_1\det(P_{S, \{1, 2\}})^2~\mbox{if}~|S|=2;\\
q_1~\mbox{if}~|S|=0;
\end{array}\right.\]
but when $k\geq 2$ then there is no clear expression for $\frac{\det(L_S)}{\det(I_n+L)}$ as a mixture of probability distributions. The problem of adapting \cite[Algorithm 4.4.2]{Hough} to the case $k\geq 2$ thus remains an open problem.

\section{Application to the simulation of attractive couplings of repulsive DPPs}\label{sec: Simulation}

\begin{figure}[h!]
\subfloat[\label{fig:example_data_1} Location of ant nests of two ant species in a site in northern Greece.]{\includegraphics[width=6.5cm]{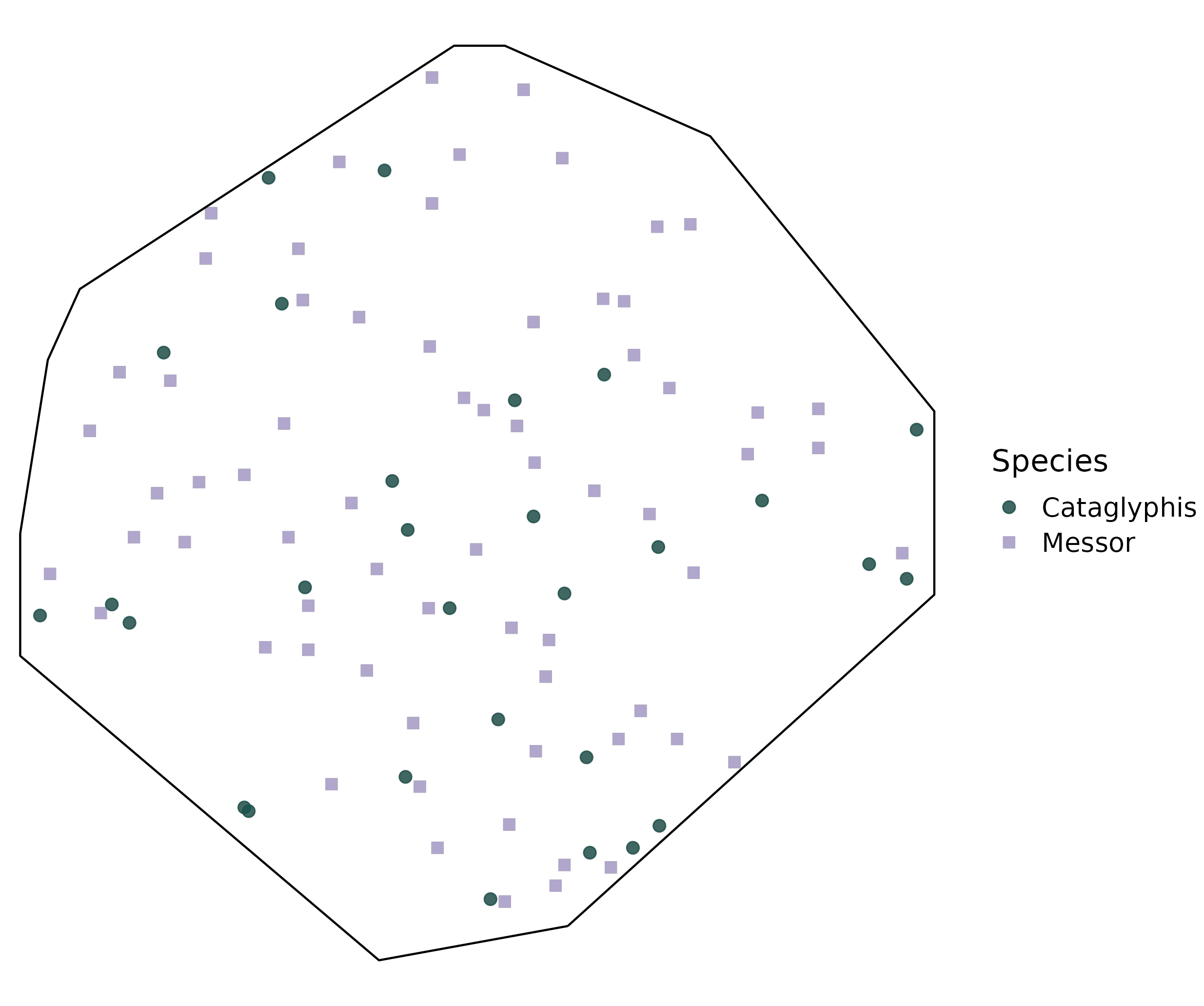}}
\hfill
\subfloat[\label{fig:example_data_2} Location of the 5G cell towers from two operators, Bouygues Telecom (circle) and Orange (rectangle), in the French city of Tourcoing. Two towers from both operators at the same location are indicated by a triangle.]{\includegraphics[width=7cm]{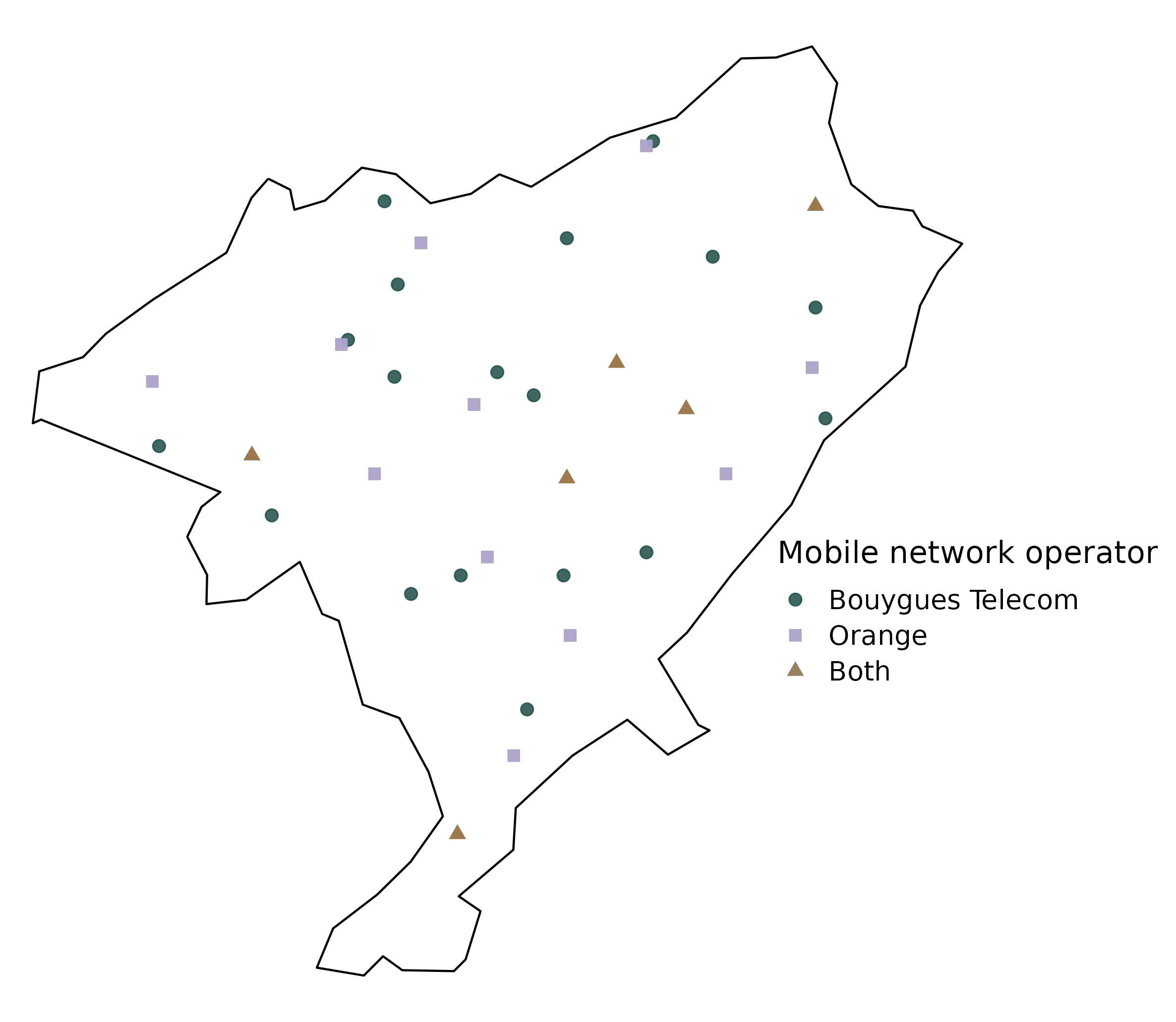}}
\caption{\label{fig:example_data} Two examples of marked point patterns with repulsion between points of the same mark and attraction between points of different marks.}
\end{figure}

We consider in this section the possible use of DPPs with nonsymmetric kernels to construct attractive coupling of repulsive DPPs. Such a distribution can be used, for example, to model marked point processes with repulsion between points of the same mark and attraction between points of different marks. One such example is the \textit{ants} dataset \cite{Harkness} in the \textit{spatstat} R package \cite{Baddeley} containing locations of nests of two species of ants, see Figure \ref{fig:example_data_1}. While ants of the same species tends to spread out their nests to avoid having to share resources, hence the repulsion in their nest location, it is suspected that one of the two species purposefully put their nest close to one of the other species since they eat other dead insects, hence the attraction between the location of nests of different species. 

\medskip

\noindent Another such example arise when looking at the location of cell towers of different network operators, see Figure \ref{fig:example_data_2}. Determinantal point processes have already been proposed to model the negative dependency in the location of base stations in a cellular network \cite{Deng, Gomez, Li}. Following a similar idea, 5G antennas owned by a given mobile network operator tends to also be well spread out spatially leading to repulsion in their locations. On the contrary, different mobile network operators tends to have antennas in similar, sometimes identical, places. This results in a positive correlation between the location of towers owned by different operators.

\noindent In both examples, DPPs with symmetric kernels are a natural candidate to model the repulsive behaviour of each mark but not the full point process. We propose instead a model for the complete marked point process using DPPs with nonsymmetric kernels.

\subsection{DPP couplings}\label{subsec: Coupling}

We begin by recalling a well-known natural way of making couplings of determinantal point processes. To express such couplings, we begin by pointing out that a coupling of two point processes on $[n]$ can be seen as a point process on $[2n]$ using the bijection
\[\begin{array}{ccc}
[2n] & \leftrightarrow & [n]\times [n]\\
X & \leftrightarrow & (X_1, X_2)
\end{array}~\mbox{where}~\left\{\begin{array}{l}
X_1 = X\cap [n];\\
X_2 = \big\{i-n, i\in X\cap \{n+1,\cdots, 2n\}\big\}.
\end{array}\right.\]
As an abuse of notation, we thus write any DPP $X$ on $[2n]$ as a coupling $(X_1, X_2)$ of DPPs on $[n]$. Now, we consider two DPP kernels $K^{(1)},K^{(2)}\in\mathcal{M}_n(\bR)$ and a $2\times 2$ block matrix of the form
\begin{equation}\label{eq: DPP_coupling}
\mathbb{K}=\begin{pmatrix} K^{(1)} & * \\
* & K^{(2)}\end{pmatrix}\in\mathcal{M}_{2n}(\bR)
\end{equation}
such that $\mathbb{K}$ is a DPP kernel. If $(X_1, X_2)\sim \DPP(\mathbb{K})$ then $X_1\sim \DPP(K^{(1)})$ and $X_2\sim \DPP(K^{(2)})$ by Proposition \ref{prop: Already_proved}~\ref{subprop: Marginal}. This way, we can write a coupling of two DPPs on $[n]$ as a DPP on $[2n]$. With a similar abuse of notation we also write $\bK_{(S, T)}$ for the principal submatrix of $\bK$ constructed from the rows and columns of $\bK$ indexed by $S\cup\{i+n, i\in T\}$.

\medskip

\noindent Such a setting has been considered before in \cite{Gourieroux} for dynamic random sets but only using symmetric kernels. This only allowed the dynamic to be repulsive since the symmetry of $\bK$ implies that
\[\mathbb{P}(i\in X_1, j\in X_2) - \mathbb{P}(i\in X_1)\mathbb{P}(j\in X_2) = -\mathbb{K}_{i, j+n}^2\leq 0\]
and therefore this coupling only allows negative cross-correlations in this case. This is why we suggest instead the possibility of using nonsymmetric matrices in order to keep this very natural way of constructing DPP couplings while allowing the coupling to have some positive cross-dependencies. While it is possible to find other ways of constructing couplings of DPPs, the fact that we can write the coupling as a DPP itself is important from a statistical standpoint since it allows the use of the nice statistical properties of DPPs like the estimation methods based on their moments \cite{Lavancier, Lavancier21}, their likelihood \cite{Poinas, Gartrell} and the various associated CLTs \cite{Biscio, Poinas17, Brunel}.

\subsection{Examples of couplings}

One of the earliest example of construction of a specific DPP coupling (using only symmetric DPP kernels) comes from \cite{Lyons} where the author considered a symmetric DPP kernel $K$ and the coupling with kernel
\[\mathbb{K}=\begin{pmatrix}K & \sqrt{K(I_n-K)}\\
\sqrt{K(I_n-K)} & I_n-K
\end{pmatrix}\]
which has the nice property of having a deterministic number of points since $\mathbb{K}$ is a projection matrix. This coupling is also an important tool to show that DPPs with symmetric kernels are strongly Rayleigh \cite{Borcea}. 

\medskip

\noindent In order to illustrate that we can use kernels \eqref{eq: DPP_coupling} to construct couplings with various degrees of repulsion and attraction we consider the following setting. Let $X$ be a DPP. We construct two new point processes $X_1$ and $X_2$ as subsets of $X$ the following way: for each point $i\in X$ let $B_i\sim b(\frac{p_i}{2-p_i})$, for some $p_i\in [0, 1]$. If $B_i=1$ then $X$ gives the point $i$ to both $X_1$ and $X_2$, otherwise $X$ gives the point $i$ to either $X_1$ with probability $1/2$ or $X_2$ with probability $1/2$. As a result, we get that $X_1, X_2\subset X$ and the higher the probabilities $p_i$ are the more points in common both $X_1$ and $X_2$ will get. We now show that $(X_1, X_2)$ can be written as a DPP coupling of the form \eqref{eq: DPP_coupling}.

\begin{proposition}\label{prop: Construction_coupling}
Let $K\in\mathcal{M}_n(\bR)$ be a DPP kernel and let $p=(p_1,\cdots p_n)\in[0, 1]^n$ such that $D(2\cdot 1_n-p)K$ is a DPP kernel. Let $X\sim \DPP(D(2\cdot 1_n-p)K)$, let $B_1,\cdots, B_n$ be Bernoulli random variables with $B_i\sim b(\frac{p_i}{2-p_i})$ and $B'_1,\cdots, B'_n$ be Bernoulli random variables with $B'_i\sim b(1/2)$. All $B_i$, $B'_i$ and $X$ are assumed to be mutually independent. We construct $X_1$ and $X_2$ the following way:
\[\begin{array}{l}
X_1=\{i\in X,~B_i=1~\mbox{or}~(B_i=0~\mbox{and}~B'_i=1)\},\\
X_2=\{i\in X,~B_i=1~\mbox{or}~(B_i=0~\mbox{and}~B'_i=0)\}.
\end{array}\]
Then, $(X_1, X_2)\sim \DPP(\bK)$ where
\[\bK=\begin{pmatrix} K & K-D(p) \\
K & K
\end{pmatrix}.\]
\end{proposition}
\begin{proof}
Let $S,T\subset [n]$. We write $A=S\cap T$ and decompose $S$ and $T$ into $S=S'\cup A$ and $T=T'\cup A$ with $A\cap S'=A\cap T'=\emptyset$. By definition of $X_1$ and $X_2$ we have
\begin{align*}
&\mathbb{P}(S\subset X_1, T\subset X_2)\\
=~&\mathbb{P}(S\cup T\subset X)\mathbb{P}(\forall i\in S',~B_i=1~\mbox{or}~(B_i=0~\mbox{and}~B'_i=1))\\
&~~~~~~~~~~~~~~~~~~~\times\mathbb{P}(\forall i\in T',~B_i=1~\mbox{or}~(B_i=0~\mbox{and}~B'_i=0))\mathbb{P}(\forall i\in A,~B_i=1)\\
=~&\det((D(2\cdot 1_n-p)K)_{S\cup T})\prod_{i\in S'\cup T'}\left(\frac{p_i}{2-p_i}+\frac{1}{2}\left(1-\frac{p_i}{2-p_i}\right)\right)\prod_{i\in A}\left(\frac{p_i}{2-p_i}\right)\\
=~&\prod_{i\in S\cup T}(2-p_i)\det(K_{S\cup T})\prod_{i\in S'\cup T'}\left(\frac{1}{2-p_i}\right)\prod_{i\in A}\left(\frac{p_i}{2-p_i}\right)\\
=~&\det(D(p)_A)\det(K_{S\cup T}).
\end{align*}
Moreover, with the right permutation of rows and columns we can write
\begin{align*}
\det(\bK_{(S, T)})&=\det\begin{pmatrix}
K_{S'} & K_{S', A} & K_{S', A} & K_{S', T} \\
K_{A, S'} & K_{A} & K_{A} - D(p)_A & K_{A, T} \\
K_{A, S'} & K_{A} & K_{A} & K_{A, T} \\
K_{T', S'} & K_{T', A} & K_{T', A} & K_{T'} \\
\end{pmatrix}\\
&=\det\begin{pmatrix}
K_{S'} & K_{S', A} & K_{S', A} & K_{S', T} \\
K_{A, S'} & K_{A} & K_{A} - D(p)_A & K_{A, T} \\
0 & 0 & D(p)_A & 0 \\
K_{T', S'} & K_{T', A} & K_{T', A} & K_{T'} \\
\end{pmatrix}\\
&=\det(D(p)_A)\det(K_{S\cup T}).\qedhere
\end{align*}
\end{proof}

\noindent Note that when $p=1_n$ in Proposition \ref{prop: Construction_coupling} then $X_1$ and $X_2$ are constructed by splitting each point of $X$ into either $X_1$ or $X_2$ with probability $\frac 1 2$. It corresponds to a very repulsive coupling since $X_1\cap X_2=\emptyset$ almost surely. This particular case matches the setting of \cite{Affandi} and we recover the following result that was pointed out in \cite{Gourieroux}:

\begin{corollary} \label{prop: Ex_exclusion}
Let $K\in\mathcal{M}_n(\bR)$ such that $2K$ is a DPP kernel. Let $X\sim \DPP(2K)$. We randomly split $X$ into two disjoint subsets $X_1$ and $X_2$ such that for each $i\in X$ we choose independently and with probability $1/2$ whether $i\in X_1$ or $i\in X_2$. Then, $(X_1, X_2)\sim \DPP(\mathbb{K})$ where
\[\mathbb{K}=\begin{pmatrix} K & K \\
K & K
\end{pmatrix}.\]
\end{corollary}

\noindent When $p=0_n$ in Proposition \ref{prop: Construction_coupling} then $X_1=X_2=X$, corresponding to the most attractive coupling possible of two DPPs, and we get the following result:

\begin{corollary}\label{prop: identical_coupling}
Let $X\sim\DPP(K)$ where $K\in\mathcal{M}_n(\bR)$ is a DPP kernel. Then $(X, X)\sim\DPP(\bK)$ where
\[\bK=\begin{pmatrix} K & K - I_n \\
K & K
\end{pmatrix}.\]
\end{corollary}

\noindent Both cases shows that the setting of \eqref{eq: DPP_coupling} allows for the construction of both very repulsive and very attractive couplings of two DPPs with the same distribution. Note that \eqref{eq: DPP_coupling} also allows the construction of an independent couplings of any DPP with kernel $K^{(1)}$ and any DPP with kernel $K^{(2)}$ since such a coupling can be written as a DPP with kernel $\mathbb{K}=\begin{pmatrix} K^{(1)} & * \\
0 & K^{(2)}\end{pmatrix}$.

\subsection{Construction of an attractive coupling of repulsive DPPs}

In order to construct an attractive coupling of two repulsive DPPs we consider $K$ to be a symmetric matrix with eigenvalues in $[0, 1]$ (and thus a DPP kernel) and we construct a coupling $(X_1, X_2)\sim \DPP(\bK)$ on $[2n]$ where $\mathbb{K}$ is of the form
\begin{equation}\label{eq: Coupling_Matrix}
\mathbb{K}=\begin{pmatrix}
K & A \\ B & K
\end{pmatrix}.
\end{equation}
We then have $X_1\sim \DPP(K)$, $X_2\sim \DPP(K)$ and since $K$ is symmetric then $X_1$ and $X_2$ are both repulsive point processes. Now, the cross dependency between $X_1$ and $X_2$ can be observed using the quantities
\[\mathbb{P}(i\in X_1, j\in X_2) - \mathbb{P}(i\in X_1)\mathbb{P}(j\in X_2) = -A_{i,j}B_{j,i}.\]
In order to have positive cross correlations we thus need to assume at minimum that $A_{i,i}B_{i,i}\leq 0$ for all $i\in [n]$. A natural way of satisfying this assumption is to take $A=-B^T$ although Proposition \ref{prop: Particle_Hole} gives
\[(X_1, X_2)\sim \DPP\left(\begin{pmatrix}
K & A \\ -A^T & K
\end{pmatrix}\right)~\Leftrightarrow~(X_1, X_2^c)\sim \DPP\left(\begin{pmatrix}
K & A \\ A^T & I_n-K
\end{pmatrix}\right).\]
Thus, the distribution of $(X_1, X_2^c)$ is a DPP with a symmetric kernel and we recover a well known setting \cite{Gourieroux}. Unfortunately, this setting is too restrictive. It does not cover the couplings presented in Proposition \ref{prop: identical_coupling} and when using it in our simulations we were only able to use such kernels to generate DPP couplings with barely any attraction. Another natural way of constructing positive cross correlations is to take $A$ (resp. $B$) to be a symmetric positive (resp. negative) matrix. By additionally assuming that $A, B$ and $K$ commute, and are therefore simultaneously diagonalizable, and then using Proposition \ref{prop: perturb_half_identity} we can choose $A$ and $B$ such that $\mathbb{K}$ is a DPP kernel with the following result.

\begin{proposition}\label{prop: Construction_Coupling}
Let $K\in\mathcal{S}^+_n(\bR)$ with eigenvalues $\lambda=(\lambda_1,\cdots,\lambda_n)\in[0, 1]^n$ and write $K=PD(\lambda)P^T$ for the spectral decomposition of $K$. Let $\mu=(\mu_1,\cdots,\mu_n)\in\bR_+^n$ and $\nu=(\nu_1,\cdots,\nu_n)\in\bR_+^n$ such that
\begin{equation}\label{eq: Bound_ev}
(\lambda_i-1/2)^2 + \frac{1}{2}\left(\mu_i^2+\nu_i^2 + |\mu_i-\nu_i|\sqrt{4\left(\lambda_i-\frac{1}{2}\right)^2+(\mu_i+\nu_i)^2}\right)\leq\frac{1}{4}
\end{equation}
and define $A=PD(\mu)P^T\in\mathcal{S}^+_n(\bR)$ and $-B=PD(\nu)P^T\in\mathcal{S}^+_n(\bR)$. Then $\mathbb{K}=\begin{pmatrix}
K & A \\ B & K
\end{pmatrix}$ is a DPP kernel.
\end{proposition}
\begin{proof}
We can write $\mathbb{K}$ as
\begin{align*}
\mathbb{K}&=\begin{pmatrix}
P & 0 \\ 0 & P
\end{pmatrix}\begin{pmatrix}
D(\lambda) & D(\mu) \\ -D(\nu) & D(\lambda)
\end{pmatrix}\begin{pmatrix}
P^T & 0 \\ 0 & P^T
\end{pmatrix}\\
&=\frac{1}{2}\left(I_{2n}-\mathbb{M}\right)~~\mbox{where}~~\mathbb{M}=-2\begin{pmatrix}
P & 0 \\ 0 & P
\end{pmatrix}\begin{pmatrix}
D(\lambda-\frac 1 2 1_n) & D(\mu) \\ -D(\nu) & D(\lambda-\frac 1 2 1_n)
\end{pmatrix}\begin{pmatrix}
P^T & 0 \\ 0 & P^T
\end{pmatrix}.
\end{align*}
We then have
\[\mathbb{M}\mathbb{M}^T=4\begin{pmatrix}
P & 0 \\ 0 & P
\end{pmatrix}\begin{pmatrix}
D((\lambda-\frac 1 2 1_n)^2+\mu^2) & D((\mu-\nu)(\lambda-\frac 1 2 1_n)) \\ D((\mu-\nu)(\lambda-\frac 1 2 1_n)) & D((\lambda-\frac 1 2 1_n)^2+\nu^2)
\end{pmatrix}\begin{pmatrix}
P^T & 0 \\ 0 & P^T
\end{pmatrix}\]
we can then deduce that the $2n$ eigenvalues of $\mathbb{M}\mathbb{M}^T$ (and thus the singular values squared of $\mathbb{M}$) are equal to four times the two eigenvalues of the $n$ matrices
\[\begin{pmatrix}
(\lambda_i-1/2)^2+\mu_i^2 & (\mu_i-\nu_i)(\lambda_i-1/2) \\ (\mu_i-\nu_i)(\lambda_i-1/2) & (\lambda_i-1/2)^2+\nu_i^2
\end{pmatrix}\]
for each $i\in[n]$. It is straightforward to show that these eigenvalues are
\[(\lambda_i-1/2)^2 + \frac{1}{2}\left(\mu_i^2+\nu_i^2\pm|\mu_i-\nu_i|\sqrt{4\left(\lambda_i-\frac{1}{2}\right)^2+(\mu_i+\nu_i)^2}\right).\]
If these quantities are all $\leq 1/4$ then Proposition \ref{prop: perturb_half_identity} shows that $\mathbb{K}$ is a DPP kernel.
\end{proof}

\noindent Note that since the first term of \eqref{eq: Bound_ev} has values in $[0, \frac 1 4]$ and the second term can be made arbitrarily small by choosing $\mu_i$ and $\nu_i$ close to $0$ then it is always possible to find some $\mu$ and $\nu$ satisfying \eqref{eq: Bound_ev}. More precisely, with some straightforward computations we can show that \eqref{eq: Bound_ev} is equivalent to
\[\max(\mu_i^2 + \nu_i^2, |\mu_i - \nu_i| + 2\mu_i\nu_i)\leq 2\lambda_i(1-\lambda_i).\]
This identity shows that the closer the eigenvalues of 
$K$ are from $0$ and $1$, the more constrained the possible eigenvalues of $A$ and $B$ will be.

\subsection{Simulation of nonsymmetric DPPs}

As mentioned in Section \ref{subsec: DPP_sim}, the main method used to simulate both discrete and continuous DPPs with symmetric kernels uses a symmetric version of Proposition \ref{prop: bad_mixing} to turn the DPP kernel into a projection matrix and then simulate projection DPPs with \cite[Algorithm 4.4.2]{Hough}. While this is the only known method used to perfectly simulate continuous DPPs, other algorithms have been developed to simulate discrete DPPs, even with generic kernels. Since we were not able to modify \cite[Algorithm 4.4.2]{Hough} for DPPs with nonsymmetric kernel we thus relied on using the generic algorithm of \cite{Poulson} for our simulations. We also mention \cite{Han, Launay} proposing different algorithms that can also be used to simulate discrete DPPs with nonsymmetric kernels.

\subsection{Numerical results}

\begin{figure}[h!]
\centering
\subfloat[Regular grid of $31\times 31$ points]{\includegraphics[width=6cm]{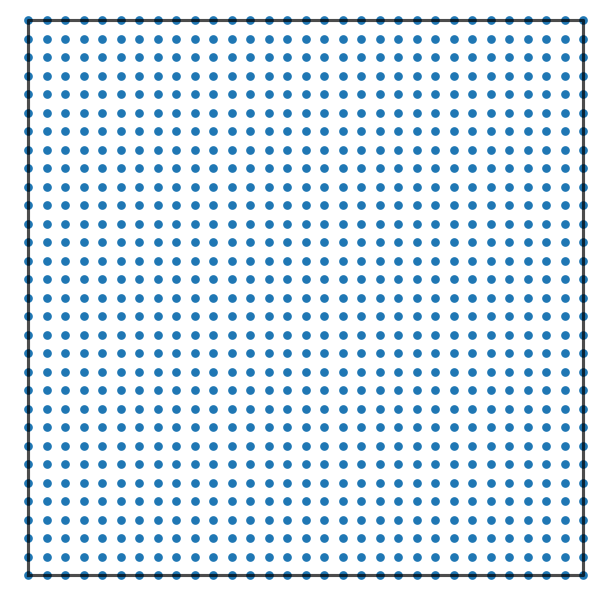}}
\hfill
\subfloat[Subset of the random grid chosen by a DPP]{\includegraphics[width=6cm]{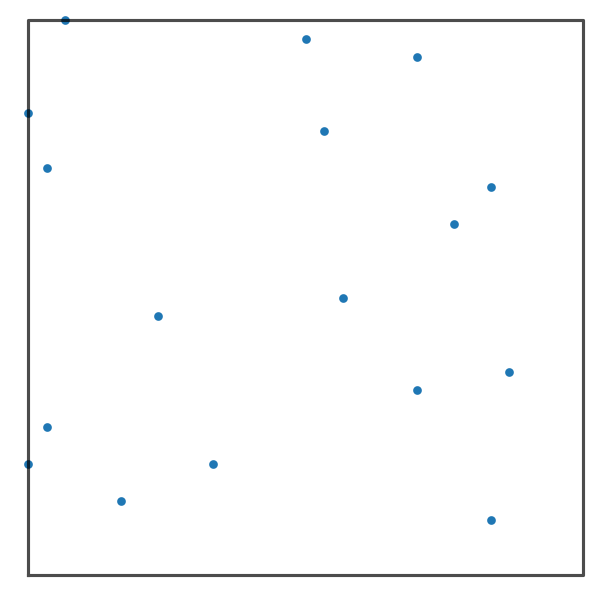}}
\caption{Example of a simulation of a DPP on a $31\times 31$ regular grid of $[0, 1]^2$ with a symmetric kernel $K$ satisfying $K_{i,j}=0.02\exp\left(-\frac{\|P_i-P_j\|_2^2}{0.018}\right).$}
\label{fig: Example_simu}
\end{figure}

For our simulations, we use DPPs to simulate random subsets of a regular grid in $[0, 1]^2$ as an approximation of a continuous DPPs on the unit box. For any $k\in\bN$ we define the regular grid of $n=(k+1)^2$ points of points of the form $(\frac{i}{k}, \frac{j}{k})\in[0, 1]^2$ and denote by $P_i$ the $i$-th point of this grid with some arbitrary ordering. We then consider kernels of the form $K_{i,j}=f(\|P_i-P_j\|)$, where $\|.\|$ is the euclidean norm on $\bR^2$ and $f$ is some well-chosen function. It was shown in \cite{Poinas} that this setting approximates the behaviour of stationary continuous DPPs. As an illustration, we give an example of the simulation of a regular DPP with a symmetric kernel in Figure \ref{fig: Example_simu}.

\begin{figure}[h!]
\subfloat[Simulation of the coupling of $X_1$ (squares) and $X_2$ (circles). \label{fig: Coupling1_simu}]{\includegraphics[height=6cm]{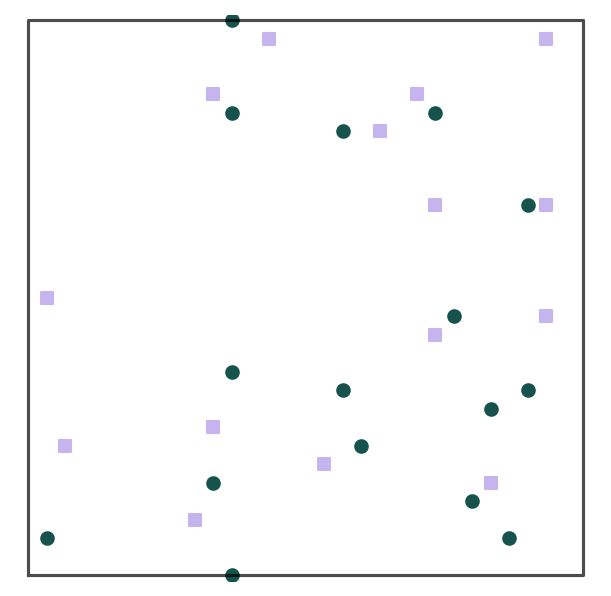}}
\hfill
\subfloat[Inclusion probability for $X_2$ given the simulation of $X_1$ obtained in Figure \ref{fig: Coupling1_simu} \label{fig: Coupling1_inclusion}]{\includegraphics[height=6cm]{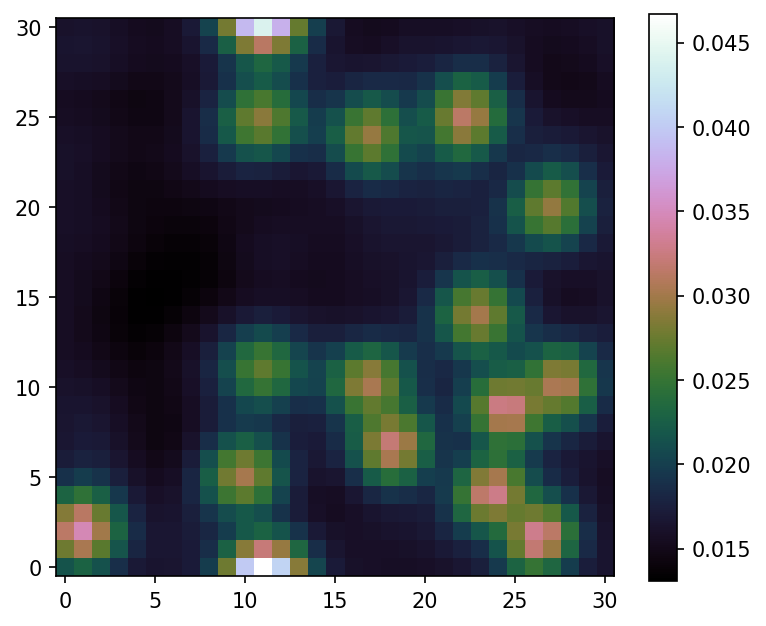}}
\caption{Illustration of an attractive coupling $(X_1, X_2)$ of two DPPs with a Gaussian Kernel.}
\label{fig: Example_simu_coupling_1}
\end{figure}

\begin{figure}[htbp]
\subfloat[Simulation of the coupling of $X_1$ (squares) and $X_2$ (circles). Points in $X_1\cap X_2$ are displayed as triangles. \label{fig: Coupling2_simu}]{\includegraphics[height=6cm]{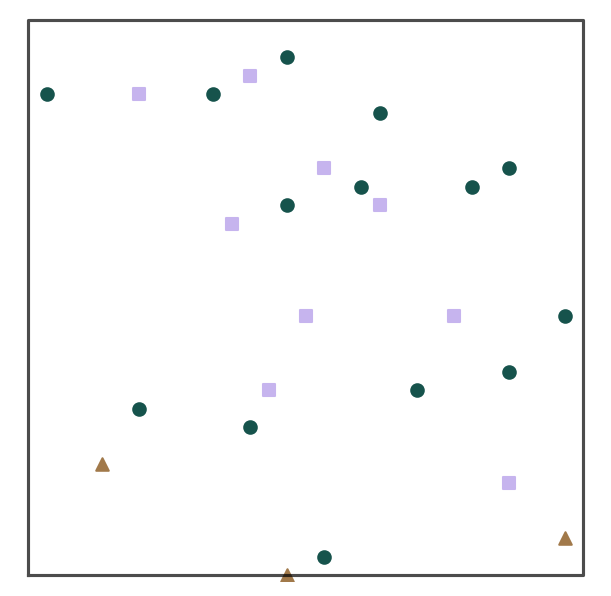}}
\hfill
\subfloat[Inclusion probability for $X_2$ given the simulation of $X_1$ obtained in Figure \ref{fig: Coupling2_simu} \label{fig: Coupling2_inclusion}]{\includegraphics[height=6cm]{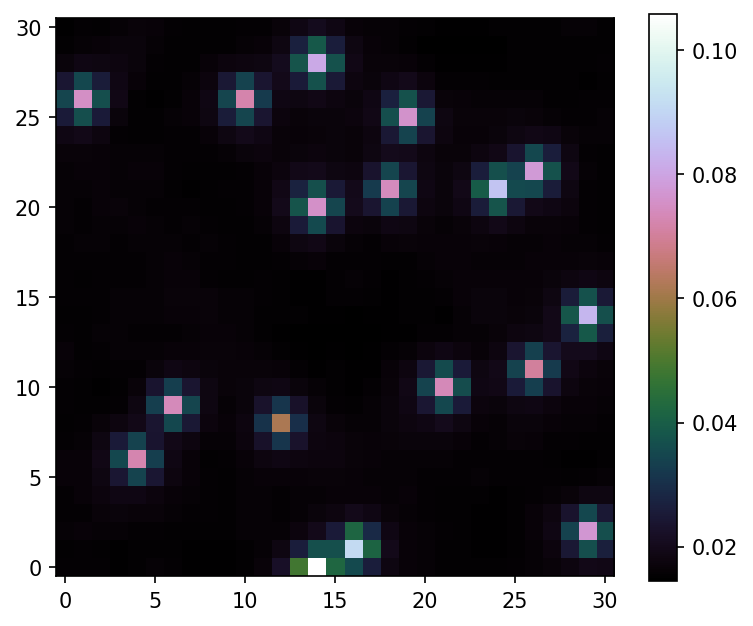}}
\caption{Illustration of an attractive coupling $(X_1, X_2)$ of two DPPs with a Cauchy Kernel.}
\label{fig: Example_simu_coupling_2}
\end{figure}

\noindent We then show in Figures \ref{fig: Example_simu_coupling_1} and \ref{fig: Example_simu_coupling_2} two attractive couplings of two repulsive DPPs with the same kernel on a regular grid of $[0, 1]^2$ with $31\times 31$ points. In both cases we simulate a coupling $(X_1, X_2)$ with kernel of the form \eqref{eq: Coupling_Matrix}. For Figure \ref{fig: Example_simu_coupling_1} we chose $K_{i,j}=0.02\exp\left(-\frac{\|P_i-P_j\|_2^2}{0.018}\right)$, called a Gaussian kernel, and for Figure \ref{fig: Example_simu_coupling_2} we chose $K_{i,j}=0.02\left(1+\left(\frac{\|y-x\|}{0.075}\right)^2 \right)^{-1.1}$, called a Cauchy kernel. In both cases, the matrices $A$ and $B$ were chosen using the method described in Proposition \ref{prop: Construction_Coupling} with all $\mu_i$ and $\nu_i$ chosen randomly among the eigenvalues satisfying \eqref{eq: Bound_ev}. Figures \ref{fig: Coupling1_simu} and \ref{fig: Coupling2_simu} both shows a simulation of the coupling $(X_1, X_2)$. In order to illustrate the attraction between $X_1$ and $X_2$ we take the simulation of $(X_1, X_2)$ then we "forget" the value of $X_2$ and we compute the probability that each point $P_i$ is in $X_2$ given the value of $X_1$ we obtained. We then present these values as a heatmap in Figures \ref{fig: Coupling1_inclusion} and \ref{fig: Coupling2_inclusion}.


\noindent Both simulations in Figure \ref{fig: Example_simu_coupling_1} and \ref{fig: Example_simu_coupling_2} display the behaviour we intended with points of $(X_1, X_2)$ being more likely to be close to each other than if $X_1$ and $X_2$ were simulated independently. Both examples also shows different kind of attraction with Figure \ref{fig: Example_simu_coupling_1} showing a small attraction with a long range, leading to points in $X_1$ and $X_2$ falling roughly in similar places as in Figure \ref{fig:example_data_1}, whereas Figure \ref{fig: Example_simu_coupling_1} show a strong attraction but with a very short range, leading to most points in $X_1$ and $X_2$ behaving like independent DPPs except for a few points that falls very close to each other or at the exact same location as in Figure \ref{fig:example_data_2}. These examples thus illustrates the possibility to use DPPs with nonsymmetric kernels in order to construct coupling with varying degrees of attraction.

\medskip

\noindent Of course, these examples are mostly a proof of concept. A lot more work is needed to see if there is a proper way to, at minimum, control the strength and range of attraction in the coupling as well as showing if this model is viable in statistical applications. But these tasks are outside the scope of this paper and they remained to be explored in future works.

\section{Discussion and open problems}

We showed in this paper the close link between the theory of $P_0$ matrices and the theory of DPPs, both area of Mathematics having some results being almost equivalent despite being expressed in a seemingly different setting. We used these close links to adapt some of the properties of $P_0$ matrices to describe the properties of DPPs with nonsymmetric kernels. We focused in particular on DPPs with kernels of the form $\frac{1}{2}(I_n - M)$ with $\|M\|_2\leq 1$ that generalize symmetric kernels, are simple to construct and for whom we adapted some of the common properties of DPPs with symmetric kernels. We finally illustrated a possible use of DPPs with generic kernels for building models of marked point processes with repulsion between points of the same mark and attraction between points of different marks.

\medskip

\noindent There is obviously still a lot more work to be done in order to deepen our understanding of the theory of DPPs so we finish by mentioning a non exhaustive list of important properties of DPPs with symmetric kernels whose potential generalization to nonsymmetric kernels would be interesting to be further investigated.

\subsubsection*{Continuous setting}

A lot of applications of DPPs, including the original purpose of modeling fermions \cite{Macchi}, consider a continuous setting, meaning that DPPs are locally finite random subsets of a general Hilbert set $H$, usually $\bR^d$. The symmetric kernel $K$ of finite DPPs is then replaced by a locally square integrable, locally of trace class hermitian operator $\mathcal{K}:H^2\rightarrow \bR$ with eigenvalues in $[0, 1]$. The immediate questions arising from the results in Sections \ref{sec: Properties} and Sections \ref{sec: prop_DPPs} is if it is possible to replace the assumption that $\mathcal{K}$ is hermitian with a continuous analogous to Proposition \ref{prop: perturb_half_identity}.

\subsubsection*{Negative association}

All DPPs with symmetric kernels satisfy a strong dependency property called negative association \cite{Borcea, Lyons} that extends to continuous DPPs with numerous applications \cite{Ghosh, Poinas17}. Since negative association implies pairwise repulsion then a necessary condition for a DPP with kernel $K$ to be negatively associated is that $K_{i,j}K_{j,i}\geq 0$ for any $i,j\in [n]$ which is obviously not satisfied by all generic DPP kernels. The question naturally arising from this conclusion is what condition on $K$ is needed for the associated DPP to be negatively associated? In fact, does a negatively associated DPP that does not have the same distribution as a DPP with symmetric kernel even exist?

\subsubsection*{Stochastic dominance}

Consider two symmetric DPP kernels $K$ and $K'$. If $K'-K$ is positive semi-definite then the DPP with kernel $K'$ stochastically dominates the DPP with kernel $K$ \cite{Lyons, Borcea}. A similar result also exists for stationary continuous DPPs \cite{Goldman}. This is an important property with various applications for continuous DPPs such as describing their reach of repulsion \cite{Bacceli} or the size of their Voronoï cells \cite{Goldman}. Now, for any generic DPP kernel $K$ and $p\in [0, 1]$ the DPP with kernel $K$ stochastically dominates the DPP with kernel $pK$ by Proposition \ref{prop: Already_proved}~\ref{subprop: Thinning} even is $K-pK$ is not positive semi-definite. One can therefore wonder if there is a simple condition on a pair of kernel $(K, K')$ that generalizes the symmetric case and indicates whether the DPP with kernel $K$ is stochastically dominated by the DPP with kernel $K'$?

%\bibliographystyle{plain}
%\bibliography{ref}

%%%%%%%%%%%%%%%%%%%%%%%%%%%%%%%%%%%%%%%%%%%%%%%%%%%%%%%%%%%%%%%%%%%
%%                                                               %%
%% Use the two commands below for producing your bibliography    %%
%% with bibtex, then comment again the commands and include the  %%
%% content of the .bbl file in this file below the commands.     %%
%%                                                               %%
%%%%%%%%%%%%%%%%%%%%%%%%%%%%%%%%%%%%%%%%%%%%%%%%%%%%%%%%%%%%%%%%%%%

%\bibliographystyle{amsplain}
%\bibliography{ref}

% add below the content of your .bbl file produced by bibtex.

%%%%%%%%%%%%%%%%%%%%%%%%%%%%%%%%%%%%%%%%%%%%%%%%%%%%%%%%%%%%%%%%%%%
%%                                                               %%
%% You may add acknowledgments (optional).                       %%
%%                                                               %%
%%%%%%%%%%%%%%%%%%%%%%%%%%%%%%%%%%%%%%%%%%%%%%%%%%%%%%%%%%%%%%%%%%%

%%%%%%%%%%%%%%%%%%%%%%%%%%%%%%%%%%%%%%%%%%%%%%%%%%%%%%%%%%%%%%%%%%%
%%                                                               %%
%% You have reached the end of your document.                    %%
%%                                                               %%
%%%%%%%%%%%%%%%%%%%%%%%%%%%%%%%%%%%%%%%%%%%%%%%%%%%%%%%%%%%%%%%%%%%

\end{document}